\documentclass[smallextended]{svjour3}       

\smartqed  

\usepackage{graphicx}          
                               
\usepackage{optidef}
\usepackage{booktabs}
\usepackage{amsmath} 
\usepackage{amsfonts}
\usepackage{color}
\usepackage{bbold}
\usepackage{tikz}
\usepackage{tikz-cd}
\usepackage{textcomp}
\usepackage{mathabx}
\usepackage{fancyvrb}

\usepackage{rotating}

\usepackage{algorithm}
\usepackage{algorithmic}

\usetikzlibrary{positioning}
\usetikzlibrary{calc}

\usepackage{float}
\newfloat{algorithm}{t}{lop}

\usepackage{threeparttable}
\usepackage{listings}
\usepackage{courier}

\graphicspath{{Figures/}}

\newcommand{\software}[1]{\texttt{#1}}

\newcommand{\R}{\mathbb{R}}

\newcommand{\Matlab}{\textsc{Matlab{}}}

\newcommand{\zeros}{0}

\DeclarePairedDelimiter\abs{\lvert}{\rvert}%
\DeclarePairedDelimiter\norm{\lVert}{\rVert}%
\makeatletter
\let\oldabs\abs
\def\abs{\@ifstar{\oldabs}{\oldabs*}}
\let\oldnorm\norm
\def\norm{\@ifstar{\oldnorm}{\oldnorm*}}
\makeatother

\newcommand{\ind}[1]{_{\mathrm{#1}}}

\newcommand{\acados}{\normalfont{\texttt{acados}}~}

\newcommand{\Jsk}{J_{\mathrm{s},k}}
\newcommand{\JsN}{J_{\mathrm{s},N}}

\begin{document}

\title{\software{acados} -- a modular open-source framework for fast embedded optimal control	
\thanks{This research was supported by the EU via ERC-HIGHWIND (259 166), FP7-ITN-TEMPO (607 957), H2020-ITN-AWESCO (642 682), by the DFG in context of the Research Unit FOR 2401, and by the  Federal  Ministry  for  Economic Affairs  and  Energy  (BMWi)  via  eco4wind  (0324125B) and DyConPV (0324166B).}
}

\titlerunning{\software{acados} -- fast embedded optimal control solvers}        
\author{Robin Verschueren
		\and Gianluca Frison
		\and Dimitris Kouzoupis
		\and Jonathan Frey
		\and Niels van Duijkeren
		\and Andrea Zanelli
		\and Branimir Novoselnik
		\and Thivaharan Albin
		\and Rien Quirynen
		\and Moritz Diehl
}


\institute{Robin Verschueren \at ABB Corporate Research, Baden-D\"attwil, Switzerland
			\and Gianluca Frison \at Dimitris Kouzoupis \at Andrea Zanelli \at IMTEK, University of Freiburg, Germany
		  	\and Niels van Duijkeren \at Department of Mechanical Engineering, KU Leuven, Belgium 
			\and Branimir Novoselnik \at Faculty of Electrical Engineering and Computing, University of Zagreb, Croatia
		  	\and Thivaharan Albin \at Institute for Dynamic Systems and Control, ETH Z\"urich, Switzerland
		  	\and Rien Quirynen \at Mitsubishi Electric Research Labs, Cambridge, MA, USA
          	\and Jonathan Frey \at Moritz Diehl \at IMTEK and Faculty of Mathematics, University of Freiburg, Germany \\ \vspace{0.3cm} \\ Note: the first four authors contributed an equal amount to the work presented in this article.
}

\date{Received: date / Accepted: date}

\maketitle

\begin{abstract}
This paper presents the \software{acados}~software package, a collection of solvers for fast embedded optimization intended for fast embedded applications.
Its interfaces to higher-level languages make it useful for quickly designing an optimization-based control algorithm by putting together different algorithmic components that can be readily connected and interchanged.
Since the core of \software{acados} is written on top of a high-performance linear algebra library, we do not sacrifice computational performance.
Thus, we aim to provide both flexibility and performance through modularity, without the need to rely on automatic code generation, which facilitates maintainability and extensibility.
The main features of \software{acados} are: efficient optimal control algorithms targeting embedded devices implemented in \software{C}, linear algebra based on the high-performance \software{BLASFEO}~\cite{Frison2018} library, user-friendly interfaces to \Matlab{} and \software{Python}, and compatibility with the modeling language of \software{CasADi}~\cite{Andersson2018}.
\software{acados} is free and open-source software released under the permissive BSD 2-Clause license.

\keywords{Direct optimal control \and Optimization algorithms}
\subclass{49-04}
\end{abstract}

\section{Introduction}

Embedded optimization, according to the definition in~\cite{Ferreau2017}, is solving optimization problems \emph{autonomously} and \emph{with limited resources}.
The topic of this article is embedded optimal control, an important class of methods within embedded optimization.
It focuses on calculating optimal decisions in order to control a dynamic system as the state changes.
There exist many algorithms for embedded optimal control and quite a few were successfully applied to real-time and embedded applications such as robotic trajectory optimization~\cite{Schulman2014}, autonomous driving~\cite{Liniger2015} and drones~\cite{Hehn2011}.
For an overview on embedded optimization, we refer the reader to~\cite{Ferreau2017}.
One of the most popular approaches in embedded optimal control nowadays is model predictive control~ (MPC) \cite{Maciejowski2002,Rawlings2017,Gruene2017}.
It is based on predicting the future behavior of a system and using this information to optimize for the action at the current time step.
In linear MPC (also called linear-quadratic MPC (LQMPC)), the constraints, including the dynamic model, are affine and the objective is quadratic.
With nonlinear MPC (NMPC), some or all of the constraints and objective are nonlinear functions\footnote{In this paper, we focus on finite-dimensional problems, i.e., the objective is not a functional.}.
A related problem is that of moving horizon estimation~(MHE) for estimating states and parameters online. We will present a unified notation of NMPC and MHE problems in Section~\ref{sec:preliminaries}.

Historically, MPC was primarily applied to systems with long timescales, most notably in chemical processing~\cite{Qin1996}, due to the fact that during each time step, a computationally costly optimization problem has to be solved.
More recent algorithmic developments~\cite{Diehl2002b,Li1989,Ohtsuka2004,Bemporad1999} and increasingly powerful embedded hardware render MPC real-time feasible for applications with shorter timescales such as autonomous driving, robotics, and avionics, as in some of the references cited above.
A big role in bringing MPC to real-time applications is played by the implementation of efficient embedded optimal control methods, giving rise to software packages such as \software{MPT}~\cite{Herceg2013} for explicit MPC, \software{qpOASES}~\cite{Ferreau2014}, an active-set solver for quadratic programming (QP), \software{FORCES}~\cite{Domahidi2012,Zanelli2017b}, an interior-point solver for quadratically constrained quadratic programming (QCQP) and nonlinear programs with optimal control structure, and the \software{ACADO Code Generation} tool~\cite{Houska2011} for tailored sequential quadratic programming (SQP) based NMPC solvers.
Other examples of nonlinear embedded optimization packages are \software{VIATOC}~\cite{Kalmari2015}, \software{GRAMPC}~\cite{Englert2019} and \software{FalcOpt}~\cite{Torrisi2018}, to which we compare \software{acados} in Section~\ref{sec:numerical_results}. A non-exhaustive list of embedded optimization software packages can be found in Table~\ref{tab:software_packages}.

Ultimately, embedded software should run on a dedicated hardware platform.
In this paper, we focus solely on central processing units (CPU) as the algorithms presented typically don't profit as much from the parallelization capabilities of massively parallel hardware platforms such as graphical processing units (GPU), tensor processing units (TPU) or field programmable gate arrays (FPGA).
Please note that we don't pin ourselves down to a specific CPU architecture: \software{acados} has been shown to work with \software{x86}, \software{x86\_64}, \software{ARMv7A}, \software{ARMv8A} and \software{PowerPC} architectures.

\begin{table}
 	\centering
	\begin{threeparttable}
	\caption{\textsc{Software packages for embedded optimization.}}
	\label{tab:software_packages}
		\begin{tabular}{|r|c|c|c|c|}
			\hline
			Software & Year & Latest Reference & Targets & License \\ \hline
			\texttt{MUSCOD-II} & 1997 & \cite{Leineweber2003} & NMPC & proprietary \\
			\texttt{HQP} & 1997 & \cite{Franke1997} & NLP & LGPL v2 \\
			\texttt{MPC Toolbox} & 1998 & \cite{MathWorks2005} & LQMPC & proprietary \\
			\texttt{AutoGenU} & 2000 &  \cite{Ohtsuka2004} & NMPC & propietary \\
			\texttt{OOQP} & 2001 & \cite{Gertz2003} & QP & propietary \\
			\texttt{MPT3} & 2003 & \cite{Herceg2013} & expl. MPC & GPL \\
			\texttt{Hybrid toolbox} & 2003 & \cite{Bemporad2005} & expl. MPC & proprietary \\
			\texttt{qpOASES} & 2006 & \cite{Ferreau2014} & QP & LGPL v2.1 \\
			\texttt{PQP} & 2008 & \cite{Cairano2013} & QP & proprietary \\
			\texttt{CVXGEN} & 2009 & \cite{Mattingley2012} & LP, QP & proprietary\\
			\texttt{ACADO Codegen} & 2009 & \cite{Houska2011} & NMPC & LGPL v3 \\
			\texttt{FiOrdOs} & 2011 & \cite{Ullmann2011} & QP & GPL v3 \\
			\texttt{muAO-MPC} & 2013 & \cite{Zometa2013} & LQMPC & 3-clause BSD \\
			\texttt{FORCES} & 2011 & \cite{Domahidi2012} & QP, QCQP & proprietary \\
			\texttt{ECOS} & 2013 & \cite{Domahidi2013a} & SOCP & GPL v3 \\
			\texttt{GRAMPC} 1.0 & 2014 & \cite{Kapernick2014} & NMPC & LGPL v3 \\
			\texttt{qpDUNES} & 2014 & \cite{Frasch2015} & LQMPC & LGPL v3 \\
			\texttt{DuQuad} & 2014 & \cite{Kvamme2014} & QP & unknown \\
			\texttt{HPMPC} & 2014 & \cite{Frison2014} & LQMPC & LGPL v2.1 \\
			\texttt{VIATOC} & 2015 & \cite{Kalmari2015} & NMPC & GPL v3 and LGPL v3 \\
			\texttt{PIPS-NLP} & 2016 & \cite{Chiang2016a} & NLP & 3-clause BSD \\
			\texttt{GPAD}\tnote{$\dagger$} & 2016 & \cite{Patrinos2014} & QP & GPL v3 \\
			\texttt{Forces NLP} & 2017 & \cite{Zanelli2017b} & NMPC & proprietary \\
			\texttt{OSQP} & 2017 & \cite{Stellato2020} & QP & Apache v2.0 \\
			\texttt{FalcOpt} & 2017 & \cite{Torrisi2018} & NMPC & MIT \\
			\texttt{HPIPM} & 2017 & \cite{Frison2020a} & LQMPC, QP & 2-clause BSD \\
			\texttt{ODYS} & 2017 & \cite{Cimini2017} & QP & proprietary \\
			\texttt{protoip} & 2017 & \cite{Khusainov2017} & NMPC & unknown \\
			\texttt{SPLIT toolbox} & 2017 & \cite{Shukla2017} & LQMPC & unknown \\
			\texttt{Control Toolbox} & 2018 & \cite{Giftthaler2018} & NMPC & 2-clause BSD \\
			\texttt{GRAMPC} 2.0 & 2018 & \cite{Englert2019} & NMPC & LGPL v3 \\
			\texttt{PANOC.jl} & 2018 & \cite{Sathya2018} & NMPC & unknown \\
			\texttt{nmpc-codegen} & 2018 & \cite{Melis2018} & NMPC & LGPL v3 \\
			\texttt{PRESAS} & 2018 & \cite{Quirynen2018a} & QP & proprietary \\
			\texttt{ParNMPC} & 2018 & \cite{Deng2018} & NMPC & unknown \\
			\texttt{qrqp} & 2018 & \cite{Andersson2018a} & QP & LGPL v3 \\
			\texttt{acados} & 2018 & \cite{Verschueren2018} & NMPC & 2-clause BSD \\
			\texttt{TMPC} & 2018 & \cite{Katliar2020} & NMPC & unknown \\
			\texttt{qpSWIFT} & 2019 & \cite{Pandala2019} & QP & unknown \\
			\texttt{OpEn} & 2019 & \cite{Sopasakis2020} & NMPC & Apache v2.0 \\
			\texttt{PolyMPC} & 2019 & \cite{Listov2020} & NMPC & MIT \\
			\texttt{QPALM} & 2020 & \cite{Hermans2020} & QP & GPL v3	\\
			\texttt{NASOQ} & 2020 & \cite{Cheshmi2020} & QP & MIT \\
			\texttt{ASIPM} & 2020 & \cite{Frey2020} & QP & proprietary \\
			\hline
		\end{tabular}
		\begin{tablenotes}
			\item [$\dagger$] As implemented by Pantelis Sopasakis.
		\end{tablenotes}
	\end{threeparttable}
\end{table}


One challenge in developing software for embedded optimal control lies in the trade-off between flexibility, memory usage and speed.
Many of the software packages mentioned in Table~\ref{tab:software_packages} are based on \emph{automatic code generation}.
One reason for that is to have self-contained efficient linear algebra routines. Often however, the size of the problem and the choice of algorithms are then fixed for one specific optimal control instance, inducing a loss of flexibility.
In some cases, a compiler or a human could generate faster or more memory-efficient code.
For example, for linear algebra operations, the recently developed high-performance linear algebra package \software{BLASFEO}~\cite{Frison2018} outperforms code-generated triple-loop linear algebra routines and state of the art \software{BLAS} and \software{LAPACK} implementations for the moderate matrix size typical for embedded optimization applications.
Since the linear algebra operations typically comprise most of the computational complexity of the algorithms considered, we can use standard compiler optimizations for the code surrounding the linear algebra routines without noticeably sacrificing performance (see Section~\ref{sec:numerical_results}).
By basing \acados on \software{BLASFEO}, we exchange self-containedness for performance and flexibility. We believe this is a better trade-off.

Another challenge for embedded optimal control software is related to the process of software development.
Often, to not sacrifice speed of execution and/or memory footprint, embedded optimal control software uses global data and suffers from tight coupling between algorithmic components.
This might lead to a codebase that is difficult to understand, maintain and extend.
We choose, as opposed to some other embedded optimal control software packages, to avoid these pitfalls by not unnecessarily sacrificing maintainability and readability of the codebase for a small gain in performance  and/or a reduction of memory footprint.
We try to achieve this goal by organizing the code in a modular fashion, with formal interfaces between the different algorithmic components, as described in Section~\ref{sec:acados}.
This allows for a straightforward way of interchanging solvers, routines, and libraries needed for the embedded control algorithm.

A final aspect of embedded optimal control software that affects flexibility, memory and runtime is the choice of modeling language and corresponding automatic differentiation tool.
Several modeling languages exist, e.g., \software{Mathematica}, \software{sympy} or the \software{MATLAB Symbolic Toolbox}.
Many of these languages make use of expression trees to represent mathematical functions, which potentially leads to a large code size, high memory usage and slow evaluation of higher-order derivatives for non-trivial models.
On the contrary, the \software{CasADi}~\cite{Andersson2018} modeling language is based on expression \emph{graphs}.
This often leads to shorter instruction sequences and to smaller, typically faster code, which makes it more suitable for embedded applications.
Also, it is free and open-source software.
For these reasons, we favor \software{CasADi} for modeling nonlinear functions and differential-algebraic equations.
Additionally, \software{acados} supports the use of hand-written or code-generated dynamic models as \software{C} source files.

In summary, the contribution of this paper is a new software package for embedded optimal control, called \software{acados}.
It offers the following main features:
\begin{itemize}
	\item efficient optimal control algorithms implemented in \software{C},
	\item modular architecture enabling rapid prototyping of solution algorithms,
	\item interfaces to \software{Python} and \Matlab,
	\item high-performance linear algebra based on \software{BLASFEO}~\cite{Frison2018},
	\item compatible with \software{CasADi} expressions~\cite{Andersson2018},
	\item deployable on a variety of embedded devices,
	\item publicly available as permissively licensed free and open-source software.
\end{itemize}

The remainder of the paper is organized as follows: in Section~\ref{sec:preliminaries}, we review some important elements of nonlinear embedded optimization algorithms relevant to \software{acados}.
We discuss recent advances in embedded optimization algorithms that motivate the development of \software{acados} in Section~\ref{sec:implementation}.
The software package \software{acados} itself is introduced in Section~\ref{sec:acados}.
Various numerical experiments including hardware-in-the-loop simulations as well as comparisons to other embedded optimization packages are presented in Section~\ref{sec:numerical_results}, and the paper is concluded in Section~\ref{sec:conclusion}.

\section{Algorithmic ingredients for embedded nonlinear optimal control}
\label{sec:preliminaries}
In this section, we first introduce the problems that typically have to be solved for NMPC and MHE.
We introduce a general formulation, that can facilitate both multiple-shooting discretized NMPC and MHE problems.
Subsequently, we describe the algorithmic incredients of SQP-type methods as implemented in \texttt{acados} for the general problem.

\subsection{Nonlinear Optimal Control}
Let us first regard the continuous-time nonlinear optimal control problem (OCP) of the form
\begin{mini!}
	{\substack{x(\cdot), z(\cdot), u(\cdot)}}
	{\int_{0}^{T} \ell(x(t), z(t), u(t)) \, \mathrm{d} t + M(x(T))}
	{\label{eq:NOCP}}
	{}
	\addConstraint{x(0) = \overline{x}_0 \label{eq:NOCP_initial_state}}
	\addConstraint{0 = f(\dot{x}(t), x(t), z(t), u(t)), \quad && t \in [0, T] \label{eq:NOCP_ODE}}
	\addConstraint{0 \geq g(x(t), z(t), u(t)), && t \in [0, T].}
\end{mini!}
In this notation, $x: \R \rightarrow \R^{n_x}$ denotes the \emph{differential states}, $z: \R \rightarrow \R^{n_z}$ are the \emph{algebraic variables} and $u: \R \rightarrow \R^{n_u}$ denotes the \emph{control inputs}.
Furthermore, we use $\ell: \R^{n_x} \times \R^{n_z} \times \R^{n_u} \rightarrow \R$ for the Lagrange term or running cost, $M: \R^{n_x} \rightarrow \R$ for the Mayer term or terminal cost.
The dynamics are modeled with a set of implicit differential-algebraic equations (DAE) with right-hand-side $f: \R^{n_x} \times \R^{n_x} \times \R^{n_z} \times \R^{n_u} \rightarrow \R^{n_x} \times \R^{n_z}$.
In the remainder of the paper, we assume the implicit DAE to be of index 1, i.e.
$\partial f / (\partial \dot{x}, \partial z)$ is invertible.
The nonlinear path constraints are given by $g: \R^{n_x} \times \R^{n_z} \times \R^{n_u} \rightarrow \R^{n_g}$, and the initial value of the states is $\overline{x}_0 \in \R^{n_x}$.
We consider the horizon length $T$ to be fixed.

\subsection{Multiple Shooting Discretization}
In \software{acados}, we discretize nonlinear OCPs with a \textit{multiple shooting} approach~\cite{Bock1984}.
We introduce a time grid $[t_0, t_1, \ldots, t_N]$ with $t_k < t_{k+1}, k=0,\ldots,N-1$, discrete state variables $x_0, \ldots, x_N$, algebraic variables $z_0, \ldots, z_{N-1}$ and controls $u_0, \ldots, u_{N-1}$.
For the control trajectory, we choose a piecewise constant control parametrization.
On each time interval $ \left[t_k, t_{k+1} \right)$, we can then write the result of the numerical simulation routine as
\begin{equation*}
\begin{bmatrix}
x_{k+1} \\ z_k
\end{bmatrix}
= \phi_k(x_k, u_k)
,\quad k=0,\ldots,N-1.
\end{equation*}
where the separate components will be denoted by $\phi_k^x(x_k, u_k)$ and $\phi_k^z(x_k, u_k)$.

The multiple shooting approach can often lead to better convergence behavior compared to single shooting, where the simulation and optimization is performed sequentially, as shown in~\cite{Bock1987}.
The resulting nonlinear programming (NLP) formulation looks as follows:
\begin{mini!}
	{\substack{x_0, \ldots, x_N, \\ z_0, \ldots, z_{N-1}, \\  u_0, \ldots, u_{N-1}}}
	{\sum_{k=0}^{N-1} (t_{k+1} - t_k) \cdot \ell(x_k, z_k, u_k) + M(x_N)}
	{\label{eq:multiple_shooting}}
	{}
	\addConstraint{x_0}{= \overline{x}_0}
	\addConstraint{\begin{bmatrix}x_{k+1} \\ z_k\end{bmatrix}}{=\phi_k(x_k, u_k), \label{eq:multiple_shooting_eq}}{\quad k=0,\ldots,N-1}
	\addConstraint{0}{\geq g_k(x_k, z_k, u_k)}{\quad k=0,\ldots,N-1}.
\end{mini!}
The optimal control formulation above is not the most general one that \software{acados} can handle: among others moving horizon estimation (MHE) problems, constraint relaxation via slack variables and equality constraints are supported as well.
The next section presents a more general class of optimization problems which is handled by \texttt{acados}.
\subsection{General Nonlinear Optimal Control structured Optimization Problem}
In order to facilitate 
many different OCP formulations occurring in practice,
\texttt{acados} uses the following general formulation of nonlinear optimal control structured optimization problems:
\begin{mini!}
	{\substack{x_0,\ldots, x_N, \\ u_0,\ldots, u_{N-1}, \\ z_0,\ldots, z_{N-1},\\ s_0, \ldots, s_N} }
	{\sum_{k=0}^{N-1} l_k(x_k, u_k, z_k)  + M(x_N) + \sum_{k=0}^{N} \rho_k(s_k)}
	{\label{eq:acados_OCP}}
	{}
	\addConstraint{\begin{bmatrix}x_{k+1} \\ z_k\end{bmatrix}\label{eq:acados_OCP_eq}}{=\phi_k(x_k, u_k),}{\quad k=0,\ldots,N-1}
	\addConstraint{0}{\geq g_k(x_k, z_k, u_k) - \Jsk s_k \label{eq:acados_OCP_ineq}}{\quad k=0,\ldots,N-1}
	\addConstraint{0}{\geq g_N(x_N) - \JsN s_N}
	\addConstraint{0}{\leq s_k}{\quad k=0,\ldots,N.}
\end{mini!}
Here, the initial state constraint from \eqref{eq:NOCP_initial_state} is contained in the inequality constraint \eqref{eq:acados_OCP_ineq} for $ k = 0 $.

Note that this general formulation also includes slack variables $ s_k $ which could alternatively be formulated as control inputs $ u_k $.
However, the slack variables $ s_k $ do not enter the system dynamics \eqref{eq:acados_OCP_eq} and can only enter the cost linearly and quadratically, i.e. via the functions $ \rho_k(\cdot) $ of the form
\begin{align*}
\rho_k(s_k) = \sum_{i=1}^{n_{s_k}} \alpha_k^i s_k^i + \beta_k^i {s_k^i}^2,
\end{align*}
with $ \alpha_k^i\in\R $, $ \beta_k^i> 0 $.
These properties motivate the separation between slack and control variables and allow for an efficient treatment in solution methods as implemented in \texttt{acados}.
The slack variables can be used to formulate soft constraints or model piecewise quadratic possibly assymetric cost functions among others.
Soft constraints are often useful in practice, for example to deal with constraint violations due to model-plant mismatch or disturbances that would yield infeasibility of \eqref{eq:multiple_shooting}.\\

Moreover, \texttt{acados} is able to handle the most common types of cost and constraint functions in a tailored fashion.
Regarding the cost functions $ l_k(\cdot) $, $ M(\cdot) $, \texttt{acados} is capable of exploiting the structure of the widely used linear and nonlinear least-squares functions, while also being able to handle general nonlinear functions.

Within the constraint functions $ g_k(\cdot) $, \texttt{acados} is able to exploit the most common constraint types, which are simply bounds on $ x_k $ and $ u_k $, linear constraints in $ x_k $ and $ u_k $, such that only truly nonlinear constraints have to be formulated and treated as such, and one could write $ g_k(\cdot) $ as
\newcommand{\Jbx}{{J\ind{bx}}}
\newcommand{\Jbu}{{J\ind{bu}}}
\begin{align}
g_k(x_k,z_k,u_k) = \begin{bmatrix} \Jbx_{,k} x_k \\ \Jbu_{,k} u_k \\ C_{\mathrm{x},k} x_k + C_{\mathrm{u},k} u_k\\ g\ind{nonl}(x_k,z_k,u_k) \end{bmatrix},
\end{align}
where $ \Jbx_{,k} $, $ \Jbu_{,k} $ are selection matrices that determine which components of $ x_k $, $ u_k $ have simple bounds.
Additionally, it is common to have upper and lower bounds on a constraint component of $ g_k $, which allows for a more efficient treatment of these constraints within \texttt{acados}.

\color{black}
The above NLP~\eqref{eq:acados_OCP} could be solved by any general-purpose NLP solver, like \software{IPOPT} \cite{Waechter2006}.
The current scope of \software{acados}, however, encompasses efficient \emph{embedded} optimal control methods for solving such structured NLPs, since these are better suited in a real-time and/or embedded setting~\cite{Diehl2009c}.
Sequential Quadratic Programming (SQP) is an example of an embeddable optimization method and will be discussed in the following section.

\subsection{Sequential Quadratic Programming and Real-Time Iterations}
We briefly present the structure of an SQP algorithm as implemented in \software{acados}.
At the very least, an embedded SQP algorithm should feature:

\begin{itemize}
	\item numerical integration of the continuous-time dynamics,
	\item generation of first-order and possibly second-order sensitivities of objective and constraints,
	\item a procedure for approximating the Hessian matrix,
	\item an efficient QP solver (typically developed separately).
\end{itemize}

Note that \emph{globalization} strategies such as line search or trust regions are considered out of scope for this paper, because globalization is not recommended in a real-time setting, since runtime bounds can not be established in general \cite{Diehl2002b}, and the initializations are typically close to the exact solution \cite{Gros2006}. In practice, globalization is typically not necessary.
Possibly extending \texttt{acados} with globalization strategies such as merit functions~\cite{Leineweber1999} and filter based line search methods \cite{Waechter2006a} for use in a fully converged setting is subject of future work.

In summary, the SQP algorithm in \software{acados} looks as follows:
\begin{align*}
w^{[i+1]} &= w^{[i]} + \Delta w^\mathrm{QP}, \qquad &&i=0,1,\ldots, \\ 
\pi^{[i+1]} &= \pi^\mathrm{QP}, \qquad &&i=0,1,\ldots, \\
\lambda^{[i+1]} &= \lambda^\mathrm{QP}, \qquad &&i=0,1,\ldots, \\
\end{align*}
where $w^{[i]}= [x_0^{[i]\top}, u_0^{[i]\top}, \ldots, u_{N-1}^{[i]\top}, x_N^{[i]\top}]^\top$ is the primal iterate at SQP iteration $i$, $\pi^{i}$ and $\lambda^{i}$ are the dual iterates, readily available from the QP solution.
Note that the algebraic variables have been eliminated from the OCP, but numerical approximations of these values are accessible from the numerical integration routine.
The step $\Delta w^\mathrm{QP}$ is computed by solving the following QP:


\begin{mini}[1]
	{\substack{\Delta x_0,\ldots, \Delta x_N, \\ \Delta u_0,\ldots, \Delta u_{N-1}, \\ s_0, \ldots, s_N}}
	{\sum_{k=0}^{N-1} \begin{bmatrix}\Delta x_k \\ \Delta u_k \end{bmatrix}^\top \overbrace{\begin{bmatrix}
	Q_k &  S_k^\top \\ S_k  & R_k\end{bmatrix}}^{H_k} \, \begin{bmatrix}\Delta x_k \\ \Delta u_k\end{bmatrix} + \begin{bmatrix}q_k \\ r_k\end{bmatrix}^\top \begin{bmatrix}\Delta x_k \\ \Delta u_k \end{bmatrix}}
	{}
	{\label{eq:OCP_QP_subproblem}}
	\breakObjective{+\Delta x_N^\top Q_N \Delta x_N + q_N^\top \Delta x_N}
	\breakObjective{+ \sum_{k=0}^{N} s_k^\top P_k s_k + p_k^\top s_k}
	\addConstraint{\Delta x_{k+1}}{=A_k \Delta x_k + B_k \Delta u_k + \bar{\phi}^x_k - x_{k+1},\quad }{k=0,\ldots,N-1}
	\addConstraint{-\bar{g}_k}{\geq G_k^x \Delta x_k + G^u_k \Delta u_k - \Jsk s_k, \quad}{k=0,\ldots,N-1}
	\addConstraint{-\bar{g}_N}{\geq G^x_N \Delta x_N - \JsN s_N}	
	\addConstraint{0}{\leq s_k}{k=0,\ldots, N.}
\end{mini}
Above, we used the shorthands $ \bar{\phi}^x_k, \, \bar{g}_k,\, k=0,\ldots, N-1$ and $\bar{g}_N$ to denote the function evaluations $\phi^x_k(x_k, u_k), \,g_k(x_k, z_k, u_k)$ and $g_N(x_N)$, respectively.
We formulate the NLP in such a way that the slack variables appear directly as tailored optimization variables.
We note that some, but not all QP solvers can deal with slack variables directly.
For those that do not, slack variables are reformulated as extra input variables.

In the remainder of this section, we discuss each of the ingredients of an efficient SQP solver, described above, separately.
Generation of sensitivities using numerical integration will be treated in Section~\ref{sec:numerical_simulation}, Hessian approximation in Section~\ref{sec:convex_hessians}, structure-exploiting QP solvers in Section~\ref{sec:qp_solvers} and real-time considerations in Section~\ref{sec:embedded_sqp}.

\subsection{Numerical simulation and sensitivities}
\label{sec:numerical_simulation}

An important part of the implementation of direct shooting methods for optimal control consists of reliably and efficiently computing numerical simulation and sensitivity results for the nonlinear system of differential-algebraic equations that represents a dynamic model for our particular system of interest.

Within the family of single-step integration methods one typically distinguishes between explicit and implicit schemes~\cite{Hairer1993}.
Well-known examples of explicit integration schemes include explicit Runge-Kutta~(RK) formulas such as explicit Euler and the RK method of order~$4$.
Explicit integration schemes are easy to implement since they rely on a direct combination of explicit evaluations of the right-hand side of the system dynamics.
Instead, implicit integration schemes result in a nonlinear system of equations that implicitly defines the numerical simulation result.
Unlike explicit integration methods, the nonlinear system in implicit integration schemes generally needs to be solved numerically using an iterative procedure such as a Newton-type method.
However, implicit formulas are very popular in practice because of their improved numerical stability properties and higher order of accuracy.
Especially in case of stiff dynamical systems and implicit or differential-algebraic equations, an implicit integration scheme should often be used~\cite{Hairer1991}.

When using these numerical integration schemes within direct multiple shooting, one additionally needs a computationally efficient and reliable way of computing first (and possibly second) order derivatives of the simulation results with respect to the state and control input values:
\begin{equation}
\dfrac{\partial \phi^x_k(x_k, u_k)}{\partial x_k}, \quad \dfrac{\partial \phi^x_k(x_k, u_k)}{\partial u_k}, \quad \sum_{i=1}^{n_x}\pi_{k,i} \dfrac{\partial^2 \phi^x_{k,i}(x_k, u_k)}{\partial^2 (x_k, u_k)},
\end{equation}
where $\pi_k \in \R^{n_x}$ is called the seed vector, for which the Lagrange multipliers are used to compute the exact Hessian of the Lagrangian.
Sensitivity propagation for direct optimal control methods is typically based on a \emph{discretize-then-differentiate} type of approach such as internal numerical differentiation~(IND) in~\cite{Bock1983}.
For the class of explicit integration methods, this concept leads to a forward or backward sensitivity propagation based on algorithmic differentiation~(AD) techniques~\cite{Griewank2000}.
In case of an implicit integration scheme, the IND approach either results in iterative differentiation techniques or a direct computation of sensitivities based on the implicit function theorem~\cite{Albersmeyer2010b}.
In addition, forward-backward propagation schemes can be derived to compute the symmetric Hessian contributions~\cite{Quirynen2017a}.

Recent work in~\cite{Quirynen2017b,Quirynen2018} proposed an algorithmic approach to embed implicit integration schemes with sensitivity analysis in Newton-type optimization for direct optimal control without the need for any iterative procedure, based on the concepts of numerical \emph{condensing} and \emph{expansion} in a lifted Newton-type optimization method~\cite{Albersmeyer2010}.

\subsection{Convex Hessian Approximation Methods}
\label{sec:convex_hessians}

\paragraph{Gauss-Newton Hessian approximation.}
In the case of a (nonlinear) least-squares objective in~\eqref{eq:acados_OCP}, e.g.
when tracking a reference, we have
\begin{align*}
l_k(x_k, u_k) &= \|r_k(x_k, u_k)\|_2^2, \quad k=0,\ldots,N-1 \\
M(x_N) &= \|r_N(x_N)\|_2^2,
\end{align*}
with $r: \R^{n_x} \times \R^{n_u} \rightarrow \R^{n_{\mathrm{r}_k}}$, $r_N: \R^{n_x} \rightarrow \R^{n_{\mathrm{r}_N}}$.
Notice that this kind of residual function is a common case in embedded optimization.

The Gauss-Newton Hessian approximation amounts to linearizing "between the norm signs".
Since no second-order sensitivities are necessary, the Gauss-Newton Hessian approximation offers a competitive alternative to SQP with exact Hessians, although it converges linearly only.
We remark that for a quadratic objective function in~\eqref{eq:acados_OCP}, the same quadratic objective arises in~\eqref{eq:OCP_QP_subproblem}, and no additional computations are needed.
For more details on Gauss-Newton methods in the context of NMPC, we refer the reader to~\cite{Gros2016}.

\paragraph{Sequential Convex Quadratic Programming (SCQP).}

A generalization to using SQP with a Gauss-Newton Hessian approximation is sequential convex quadratic programming~\cite{Verschueren2016}.
In a sense, a Gauss-Newton SQP algorithm neglects any curvature present in the inequality constraints by linearizing them.
In practice however, convex-over-nonlinear objectives and/or constraints arise often, which are of the form $\varphi(c(x, u))$ with a convex function $\varphi(\cdot)$ and a nonlinear function $c(\cdot)$.
Examples include ellipsoidal terminal constraints to ensure stability of an NMPC scheme, the friction ellipse in automotive applications, or tunnel-following for robotic manipulators.

In SCQP, we still linearize the inequalities, but bring the convex contributions from the inequality constraints, multiplied with a Lagrange multiplier, into the Hessian approximation. Doing so, the SCQP Hessian is guaranteed to be positive semi-definite.
For problems that feature convex-over-nonlinear constraints, this Hessian contribution offers better convergence guarantees than a Gauss-Newton Hessian~\cite{Verschueren2016}.

\paragraph{Structure-preserving convexification with minimal regularization.}

In the last two paragraphs, we devised two Hessian approximations which are convex \emph{by construction}.
However, when the exact Hessian of the Lagrangian is used, it might be indefinite.
When this happens, the optimal direction $\Delta w_\mathrm{QP}$ cannot be guaranteed to be a descent direction.
Furthermore, many QP solver codes expect a positive (semi-)definite Hessian, even if the second order sufficient conditions for optimality are met.
The aim of regularization is to obtain an approximation $\widetilde{H} = \mathrm{blkdiag}(\widetilde{H}_0,\ldots,\widetilde{H}_N)$ with each $\widetilde{H}_k \succ 0$.
We very briefly discuss three different methods here and compare their convergence in an SQP-type setting in Section~\ref{sec:numerical_results}.

Let $V_k D_k V_k^{-1}$ be the eigenvalue decomposition of $H_k$, for $k=0,\ldots,N$.
Two simple ways of regularizing the Hessian blocks are
\newcommand{\maxdiag}{\mathrm{maxdiag}}
\begin{subequations}
	\label{eq:regularization}
	\begin{alignat}{1}
	\widetilde{H}_k &= \mathrm{project}(H_k,\epsilon) := V_k \, \big[\maxdiag(\epsilon, D_k) \big] \, V_k^{-1}, \\
	\widetilde{H}_k &= \mathrm{mirror}(H_k,\epsilon) := V_k \, \big[\maxdiag(\epsilon, \mathrm{abs}(D_k)) \big] \, V_k^{-1},
	\end{alignat}
\end{subequations}
with $\epsilon$ small, $ \mathrm{abs}(\cdot)$ defined elementwise and $ \maxdiag(\cdot, \cdot) $ selecting the maximum values on the diagonal while not changing the off-diagonal elements.
\\
Another approach is to exploit the optimal control problem structure of~\eqref{eq:OCP_QP_subproblem}.
One approach of this kind is called \emph{convexification} and has been proposed in~\cite{Verschueren2017}.
The difference with more naive regularization methods as stated above, is that it first exploits as much convexity as possible through the optimal control structure, before regularizing the remaining negative directions.
The complexity of the regularization method is linear in the horizon length, and under some conditions, the SQP iterates converge quadratically to a local solution.
An efficient implementation of this new Hessian regularization method is included in \software{acados}.
A numerical example is given in Section~\ref{sec:numerical_results}, which shows the superior convergence behavior of convexification with respect to the mirror and project regularization techniques.

\paragraph{Further convex Hessian approximations} Other Hessian approximations exist in the literature, such as BFGS and its limited-memory variant (see \cite{Nocedal2006}). They might be added to \software{acados} in the future.

\subsection{Structure-exploiting embedded QP solvers}
\label{sec:qp_solvers}

There exist different solution strategies for QP~\eqref{eq:OCP_QP_subproblem}, which we briefly describe in this section.
We note that linear-quadratic optimal control problems can be efficiently solved with the same kind of QP solving strategies presented.
As such, \software{acados}, conceived as a modular software package, can also be used to facilitate solving linear-quadratic QPs, arising in linear-quadratic MPC.\\
In the following subsections, we give an overview on sparsity exploitation for optimal control structred QPs.

\paragraph{Sparse approach.} OCP~\eqref{eq:OCP_QP_subproblem} can be solved directly by using a general-purpose sparse QP solver, e.g., \software{CVXGEN}~\cite{Mattingley2012}, \software{OOQP}~\cite{Gertz2003}, both primal-dual interior-point solvers, or \software{OSQP}~\cite{Stellato2020}, a first-order method. First-order methods are mainly based on either the fast gradient method or the alternating direction method of multipliers (ADMM).
The strict real-time requirements for solution methods make first-order methods a viable candidate.
However, they might suffer from slow convergence rates.
For an overview of first-order methods in the context of embedded optimal control, see~\cite{Ferreau2017,Kouzoupis2015}.

\paragraph{Structured approach.} OCP~\eqref{eq:OCP_QP_subproblem} is solved by exploiting its multi-stage structure, but dense linear algebra is used.
An example is the approach from~\cite{Steinbach1995,Rao1998} and in solvers like \software{FORCES}~\cite{Domahidi2012,Zanelli2017b}, \software{HPMPC}~\cite{Frison2014} and \software{HPIPM}~\cite{Frison2020a}, which are all interior-point solvers.

\paragraph{Condensing approach.} An alternative to the previous approaches is the so-called \emph{condensing} approach~\cite{Bock1984}.
By eliminating the state variables by means of the dynamic equality constraints in~\eqref{eq:OCP_QP_subproblem}, we obtain a smaller QP with only the control inputs and possibly the initial state as optimization variables.
Any general-purpose dense QP solver can then be used to solve the smaller QP, e.g. an active-set solver like \software{qpOASES}~\cite{Ferreau2014} or an interior-point method like the dense variant in \software{HPIPM}.
Since \texttt{qpOASES} is able to reuse information (i.e. \emph{warm-starting}) from one problem to the next, it is particularly well-suited for (N)MPC.
Condensing is shown to be of quadratic complexity in the horizon length~\cite{Frison2015a}.

\paragraph{Partial condensing.} A mix between the two previous approaches (i.e.
structured, condensing) for solving~\eqref{eq:OCP_QP_subproblem} can be obtained by not eliminating all state variables, but only per blocks of $N / N_2$ stages (we assume for simplicity that $N$ is an integer multiple of $N_2$), where $N_2$ is the `new' horizon length of the partially condensed problem.
By this additional degree of freedom, partial condensing enables us to find a better trade-off between horizon length and number of optimization variables, for a given problem.
For more details on partial condensing, the reader is referred to~\cite{Axehill2015}.

\paragraph{Further condensing strategies.}
There are further condensing strategies in the literature, which might be added to \acados in the future.
A method called \textit{complementary condensing} was proposed in \cite{Kirches2012b}, where the KKT system is solved in the space of the multipliers corresponding to the equality constraints.
This approach favorable if the problem has more control inputs than states.

\paragraph{Dual-Newton strategy.}

The dual Newton strategy is an algorithm that is based on dual decomposition tailored to linear-quadratic OCPs in the form of~\eqref{eq:OCP_QP_subproblem}, with an open-source implementation in the software \software{qpDUNES}~\cite{Frasch2015}.

A main advantage of the dual Newton strategy, as in most active set methods, is warm-starting. Contrary to \software{qpOASES}, \software{qpDUNES} can perform multiple active set changes per iteration.
However, a premature termination of the algorithm does not return a meaningful solution as in the case of~\software{qpOASES} since it is not feasible nor optimal for a neighboring problem.
As recently observed in~\cite{Kouzoupis2015b}, the convergence of~\software{qpDUNES} can benefit significantly from a partial condensing preprocessing step.

\subsection{Real-time iterations}
\label{sec:embedded_sqp}

In a real-time control setting, we solve NLP~\eqref{eq:acados_OCP} in sequence and under stringent time conditions.
Since the environment is anyway changing continuously, it is often sufficient to solve it approximately -- it is of no use to have a high-accuracy but past-the-deadline solution.

One such online method is the real-time iteration (RTI) scheme~\cite{Diehl2002b}.
It solves an inequality-constrained QP in each iteration.
The resulting generalized predictor is better suited for predictions across active set changes, than e.g. a tangential predictor obtained from an interior-point method.
For a brief overview, we refer the reader to~\cite{Diehl2009c}.

In each RTI, one full iteration of an SQP-type scheme is performed, including generation of the sensitivities w.r.t.
all variables.
We could introduce additional approximations by not updating all sensitivities in each RTI~\cite{Zanelli2019}.
Such approximations exist on different levels: from only updating the initial state constraint with the current estimate of the state of the system, over updating the right-hand side of the (in)equality constraints, to the full RTI.
By interleaving different approximations at different sample times, we obtain a \emph{multi-level} iterations scheme, as introduced by~\cite{Bock2005}.

\section{Algorithm implementations in \texttt{acados}}
\label{sec:implementation}

In this section, we focus on the algorithm implementations in \software{acados}.
Since \software{acados} builds on other software to handle basic linear algebra operations (\software{BLASFEO}~\cite{Frison2018}) and QPs (\software{HPIPM}~\cite{Frison2020a}), these will be presented first.
Afterwards, a short description of integrators and SQP-type optimization solvers will be given.

\subsection{Linear algebra: \software{BLASFEO}}

At the heart of all embedded optimization routines lies either an implementation of a small set of linear algebra routines (e.g.
matrix-matrix multiplication, Cholesky decomposition), or a call to a specialized linear algebra library (e.g.
\software{BLAS} and \software{LAPACK}).
Generally \software{BLAS} implementations focus on performance for large dense matrices, as used in high-performance computing and data science.
Considerably less investigated are \software{BLAS} and \software{LAPACK} implementations for small dense matrices.

Often, the linear algebra code in embedded optimization packages, for example in the~\software{ACADO Code Generation} tool or \software{Forces Pro}, is code-generated.
For very small matrix sizes (e.g.
$4 \times 4$), this technique outperforms optimized linear algebra libraries.
Furthermore, code-generation has the advantage that the code can be kept `library-free'.
However, for larger matrix sizes (e.g.
in the range $10 \times 10$ to $100 \times 100$, typical in MPC applications), code-generated linear algebra routines underperform with respect to optimized libraries.

\software{BLASFEO}~\cite{Frison2018} is a linear algebra package that targets computations for small matrices.
It offers highly optimized linear algebra routines (e.g.
\software{dgemm}, \software{dsyrk}, \software{dpotrf}), tailored for the matrix sizes typically encountered in embedded optimization.
These routines exploit architecture-specific vector instructions for floating point operations (e.g.
\software{AVX}), and focus on performance for matrices fitting in cache.

Furthermore, \software{BLASFEO} defines a packed matrix format (called panel major) which optimizes the cache usage, guaranteeing close to peak performance for matrices of sizes up to a couple hundreds.
All high-performance \software{BLASFEO} routines use this panel major matrix format, and there is a rich set of auxiliary routines to operate on this matrix format, as well as to convert from/to column- or row-major formats.
In this sense, \software{BLASFEO} provides a complete linear algebra framework, which can be used to implement many fast optimization algorithms.

Except for trivially small matrices, \software{BLASFEO} enables a considerable speedup (up to $10\times$ for some matrix sizes) in the matrix computations, compared to code-generated linear algebra kernels.
For all small matrix sizes up to, say, $300 \times 300$, \software{BLASFEO} offers a considerable speedup compared to state-of-the-art \software{BLAS} implementations, like \software{OpenBLAS}, too.
The use of \software{BLASFEO} is one of the factors why \software{acados} performs better than \software{ACADO} on medium-scale problems, as we will see in Section~\ref{sec:numerical_results}.

\subsection{Quadratic programming: \software{HPIPM}}

In the implementation of SQP-type algorithms, QP sub-problems need to be solved efficiently at each iteration.
The QP sub-problem solution is typically one of the two most expensive steps in SQP schemes, the other being the simulation and sensitivity computations of dynamical systems.

\software{HPIPM}~\cite{Frison2020a} is a library defining three QP types (dense QP, OCP QP and tree-structured OCP QP, all handling soft constraints via tailored slack variables), and a rich set of routines to create, manage and solve the QPs.
All QP solvers are Mehrotra's type primal-dual interior point methods, and they are implemented using the \software{BLASFEO} linear algebra framework.

Furthermore, \software{HPIPM} provides a set of efficiently implemented routines to convert between the different QP types, e.g. condensing routines convert an OCP QP into a dense QP, partial condensing routines convert an OCP QP into another OCP QP with shorter horizon length.
Additional expansion routines recover the original QP solution from the (partially) condensed one.

In \software{acados}, the QP framework is based on \software{HPIPM}, in the sense that \software{HPIPM} provides both the dense QP and OCP QP definitions, as well as (partial) condensing algorithms to convert them and interior-point methods to solve them.
Numerous other QP solvers are then interfaced to \software{acados} to alternatively solve the same types of QP problems.
At the time of writing, additional interfaces exist to~\software{HPMPC}, \software{qpDUNES}, \software{qpOASES}, and \software{OSQP}.

\subsection{Algorithmic Differentiation}

The main difference between software packages for linear quadratic MPC and nonlinear MPC is that NMPC packages need a modeling language for the nonlinear functions. We support \software{CasADi}, a graph-based source transformation algorithmic differentiation tool~\cite{Andersson2018}. In our workflow, a user would typically specify the dynamic continuous-time models and nonlinear constraint functions with \software{CasADi} in a high-level language such as \Matlab{} or \software{Python}.
In combination with \software{CasADi}'s code generator we obtain fast, embeddable code, see Section \ref{sec:high_level_interfaces}.

\subsection{Numerical simulation}

\software{acados} features different kinds of numerical simulation routines.
There are implementations of explicit and implicit Runge-Kutta integrators available, both of which support the optional propagation of first-order forward and adjoint sensitivities, as well as second-order sensitivities.
The explicit integrators can be used with explicit ODE models and supports different Butcher tableaus, including Euler's method and RK4.
Moreover, the implicit integrators can be used with an index-1 differential-algebraic equation (DAE) or implicit ODE model and use the Gauss-Legendre Butcher tableaus.
A novel implementation for lifted collocation integrators~\cite{Quirynen2017b} has been made part of \software{acados}, as well as a recently proposed structure-exploiting IRK algorithm, the so-called GNSF-IRK scheme~\cite{Frey2019}, discussed next.

The concept of GNSF-IRK is to rigorously exploit the linear dependencies within the dynamic system.
It extends the ideas of the linear input and linear output subsystems that have been implemented within the \software{ACADO Code Generation} tool~\cite{Quirynen2013a} and uses a more flexible structured dynamic system formulation that can also handle index-1 DAEs.
A main challenge for structure-exploiting integrators is to appropriately reformulate the dynamic system of interest into the desired structured form.
\software{acados} features an automatic transcription method for the GNSF structure~\cite{Frey2019}, implemented as a \Matlab~function for \software{CasADi} models.

A last important feature of \software{acados} is that integrators can vary from stage to stage, with e.g.
different state and control dimensions, different integration step length, or different integration schemes.

We note that all integrator modules, except for the ERK integrator, are based on hardware tailored linear algebra routines in \software{BLASFEO} to speed up the LU factorizations and the corresponding triangular system solutions, as discussed in~\cite{Frison2017}.

\subsection{SQP-type methods}

For nonlinear programming, \software{acados} offers different SQP-like methods.
A full-step SQP method is available, with different algorithmic options.
As Hessian approximations, we have Gauss-Newton Hessians, SCQP, and exact Hessians with regularization/convexification as discussed in Section~\ref{sec:embedded_sqp}, and we allow for user-defined Hessian approximations.
For use in an online setting, e.g.
in NMPC, a specialized RTI routine is available.

Both the SCQP algorithm and the convexification method of Section~\ref{sec:convex_hessians} are novel features, to the authors' knowledge, not present in any other NMPC software packages.

\section{The \software{acados} software package} \label{sec:acados}

\software{acados} implements some of the optimization methods mentioned in the previous sections.
\software{acados} is meant to be user-friendly at a high level, and efficient at a low level.
In order to balance these properties, we developed a core library written in \software{C} which exposes functionality to the \software{Python} and \Matlab{} interfaces.
In this section, we first discuss the functionality of this inner core module, we then describe internal and external interfaces that are crucial for usability.

\subsection{The \software{acados} core library}

The embedded optimization algorithms discussed in Section~\ref{sec:preliminaries} are implemented in \software{acados} in a modular fashion.
For example, there is a clear interface between an NLP solver and an integrator.
The integrator expects a linearization point $w^{[i]}$ and returns the end state of a simulated trajectory, and optionally first- and second-order sensitivities:

~

\begin{tikzpicture}
\begin{tikzcd}[every arrow/.append style={shift left}]
	\mathrm{NLP~solver} \arrow{r}{\mathrm{lin.~point}} &\mathrm{integrator} \arrow{l}{\mathrm{sim., sens.}}
\end{tikzcd}
\end{tikzpicture}

Similar diagrams can be drawn for all other algorithmic components, including (partial) condensing, QP solvers, function evaluations etc.
Each of these algorithmic components are modeled within \software{acados} as separate \emph{modules}.
Some modules can be used as standalone modules, or in combination with others.
For instance, depending on the choice of algorithm, an NLP solver will make use of some or all of the other modules.
In Table~\ref{tab:modules}, we see an overview of all modules currently present in \software{acados}, together with the implemented algorithmic variants.

\begin{table}
	\caption{\textsc{Overview of the software modules present in \software{acados}.}}
	\label{tab:modules}
	\centering
	\begin{tabular}{|c|c|}
		\hline 
		Module & Variants \\ 
		\hline \hline
		OCP QP & \begin{tabular}{@{}c@{}}
			\software{HPIPM} \\ \software{qpDUNES} \\ \software{HPMPC} \\ \software{OSQP} 
		\end{tabular} \\ 
		\hline
		Dense QP & \begin{tabular}{@{}c@{}}
			\software{HPIPM} \\ \software{qpOASES} 
		\end{tabular} \\
		\hline
		Condensing & \begin{tabular}{@{}c@{}} Full condensing (\software{HPIPM}) \\ Partial condensing (\software{HPIPM}) \end{tabular} \\
		\hline
		Simulation & \begin{tabular}{@{}c@{}} ERK \\ IRK \\ GNSF-IRK \\ lifted IRK \end{tabular} \\ 
		\hline 
		OCP NLP & \begin{tabular}{@{}c@{}} Gauss-Newton SQP \\ Gauss-Newton SCQP \\ Exact-Hessian SQP \\ RTI \end{tabular} \\ 
		\hline
		Regularization & \begin{tabular}{@{}c@{}} Projection \\ Mirroring \\ Convexification \end{tabular} \\
		\hline
		Nonlinear function & \begin{tabular}{@{}c@{}} \software{CasADi} generated functions \\ \software{C}-code functions \end{tabular} \\
		\hline
	\end{tabular} 
\end{table}

It is an important design choice that all modules are identical in their signature.
That way, all modules look similar to the users of \texttt{acados}.
For developers, it should be straightforward to extend \texttt{acados} with another module.
The signature is as follows (in \texttt{C} syntax):

\begin{lstlisting}[basicstyle=\footnotesize\ttfamily,language=C++,emph={void, int},emphstyle={\color{blue}}]
int <solver>(void *config,
             void *dims,
             <module>_in *in,
             <module>_out *out,
             void *opts,
             void *mem,
             void *work);
\end{lstlisting}

Here, \texttt{<module>} stands for the name of the module at hand, for example \texttt{ocp\_qp} for QP problems with optimal control structure or \texttt{sim} for integration problems, and \texttt{<solver>} is a placeholder for a function implementing the specific solver for problems corresponding to this module, e.g., \texttt{ocp\_qp\_hpipm} (interface to \software{HPIPM} solver) or \texttt{sim\_erk} (explicit Runge-Kutta method), etc.
Each module returns an \texttt{int} which denotes a solver-specific error status -- zero means successful completion by convention.
All input arguments are pointers.
Each of the arguments comes with a set of helper functions, called \texttt{...\_calculate\_size}, computing the size (in bytes) of the \texttt{struct} pointed to, as well as a set of functions, called \texttt{...\_assign}, to initialize a block of memory.

Some modules comprise other modules.
For example, an SQP solver for optimal control problems might need an integrator, which is on its own a proper \software{acados} module.
In this context, we call the integrator a \emph{submodule}.
Each of the arguments above, \texttt{dims, in}, etc., have fields corresponding to submodules.
As an example, the relation between an NLP solver and its submodules is depicted in Figure~\ref{fig:submodules}.
We remark that the calculation of the memory size of a module with submodules is done recursively, i.e., calling the \texttt{calculate\_size} function on the top module returns the required memory size of the top module and all of its submodules, and submodules of submodules, etc.
This allows users to allocate all the memory outside of \software{acados}, by design.

\tikzstyle{module}=[rectangle,draw=black!50,fill=gray!20,node distance=18pt]
\begin{figure}
	\centering
	\begin{tikzpicture}
	\node[module] (dyn) [align=center] {dynamics \\ (continuous)};
	\node[module] (int) [below=of dyn, align=center] {simulation \\ (ERK)};
	\node[module] (ext1) [below=of int, align=center] {external \\ function};
	\node[module] (cos) [right=of dyn, align=center] {cost \\ (nonlinear)};
	\node[module] (ext2) [below=of cos, align=center] {external \\ function};
	\node[module] (con) [right=of cos, align=center] {constraints \\ (nonlinear)};
	\node[module] (ext3) [below=of con,align=center] {external \\function};
	\node[module] (qps) [right=of con, align=center] {OCP QP solver \\ (HPIPM)};
	\node (mid) at ($(cos)!0.5!(con)$) {};
	\node[module] (NLP) [above=of mid, align=center] {OCP NLP solver \\ (SQP)};
	\draw -> (dyn) -- (NLP);
	\draw -> (cos) -- (NLP);
	\draw -> (con) -- (NLP);
	\draw -> (qps) -- (NLP);
	\draw -> (int) -- (dyn);
	\draw -> (ext1) -- (int);
	\draw -> (ext2) -- (cos);
	\draw -> (ext3) -- (con);
	\end{tikzpicture}
	\caption{Example of the relation between modules and submodules in \software{acados} for a specific case of a possible SQP algorithm.}
	\label{fig:submodules}
\end{figure}
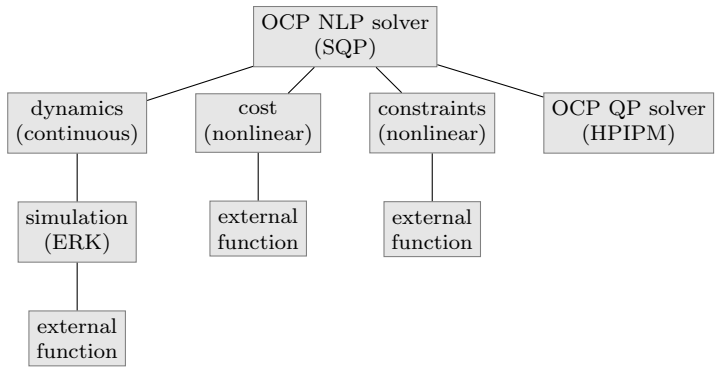

The \emph{core} library of \texttt{acados} contains mostly what has been described in this section: a collection of modules, each with corresponding data types and variants of solvers, as well as helper functions for memory management.
Using the core library directly can be cumbersome and error-prone, as many details need to be taken into account: it is designed to be efficient and flexible.
To cater to the specific needs of the end user, we offer different interfaces to the core of \software{acados}, which are described next.

\subsection{The C interface}

The \texttt{C} interface is responsible for encapsulating the low-level constructs of the \software{acados} core.

\paragraph{Choosing Solvers.} When working with the core library, all functions are specific to one variant of a module: when solving a QP with, say, \texttt{qpOASES}, the code will refer to \texttt{struct}s like \texttt{dense\_qp\_qpoases\_memory}, \texttt{dense\_qp\_qpoases\_opts}, etc.
We provide an abstraction layer to facilitate switching solvers easily.
To this end, for each module we define a `plan'.
A plan is a \texttt{struct} that contains a number of fields representing the choice of a particular combination of solvers.
For example, the plan for an SQP-type method with Gauss-Newton Hessian approximation, for a problem discretized with an ERK integrator using \texttt{HPIPM} as an underlying QP solver, reads as
\begin{lstlisting}[basicstyle=\footnotesize\ttfamily,language=C++,emph={void, int},emphstyle={\color{blue}}]
ocp_nlp_solver_plan plan = {
    {PARTIAL_CONDENSING_HPIPM},
    {ERK, ERK, ERK, ...},
    SQP_GN,
    {LINEAR_LS, ...},
    {CONTINUOUS_MODEL, ...},
    {BGH, ...},
};
\end{lstlisting}

Here, the arrays should be of the correct length (omitted for brevity, with a slight abuse of notation).
As a general rule, solvers that make use of other modules should include them in their plan.

\paragraph{Passing options.}
For each module, we manipulate a specific options \texttt{struct} using functions that take a textual representation of the option via a string.
The string encodes both the module that the option belongs to, as well as the name of the option.

\paragraph{Memory management.} Allocating memory `manually' as described above can quickly become cumbersome.
For this reason, we make available a few routines that automate that process.
To this end, in the \texttt{C} interface each module from the core library is mirrored by an additional function with signature
\begin{lstlisting}[basicstyle=\footnotesize\ttfamily,language=C++,emph={void, int},emphstyle={\color{blue}}]
<module>_solve(<module>_solver *solver,
               <module>_in *in,
               <module>_out *out);
\end{lstlisting}

The argument \texttt{solver} is a pointer to a \software{C} structure that encapsulates the data needed other than input and output.
By doing so, we reduce the amount of boilerplate code.
%
%
%
%


\paragraph{Convenience routines.}

The \software{C} interface additionally offers helper routines, so-called `setters' and `getters', that wrap the handling of the low-level \texttt{struct}s of the \software{acados} core.

\subsection{High-level interfaces}
\label{sec:high_level_interfaces}
Often, software for NMPC is coded in scripting languages.
Therefore, we offer interfaces to two popular languages for scientific computing: \software{Python} and \Matlab, where the interface to \Matlab{} is largely compatible with its free and open-source alternative \software{Octave}.
As such, we created a small domain-specific language within each of these frameworks.
We build on top of code from the \texttt{C} interface of \software{acados}.


In order to formulate the OCP \eqref{eq:acados_OCP} through the \texttt{acados} modules (cost, constraints and dynamics), the main challenge is to pass the generally nonlinear functions and their derivatives to these modules. 
The \software{Python} and \Matlab{} interfaces of \texttt{acados} use \texttt{CasADi} as a modeling language, i.e. to formulate all generally nonlinear parts of the OCP.
The \software{acados} high level interfaces are able to use \texttt{CasADi}'s code generation and algorithmic differentiation to generate the \texttt{C} functions needed for each \software{acados} module.
A major benefit of using \texttt{CasADi} as a modeling language is that the solution behavior of \texttt{acados} can be easily compared with the solutions coming from the numerous optimization tools interfaced with \texttt{CasADi}.

Once the OCP to be solved is described through the domain-specific language implemented by the high-level interfaces, a human readable self-con{\-}tained \texttt{C} project that makes use of templated code can be generated.
The generated project contains all the \texttt{C} code necessary for function and derivative evaluations generated through \texttt{CasADi} and the \texttt{C} code necessary to set up the NLP solver using the acados \texttt{C} interface.
Moreover, a \Matlab{} S-Function and a build system for its compilation is generated.
Note that this kind of code generation is inherently different from the one in \software{ACADO}, since the templated code uses only the functions exposed by the \texttt{C} interface of \texttt{acados}.
In contrast to this, \software{ACADO} generated solvers are standalone \software{C} projects that are extremely problem specific and do not rely on a common library.

With the workflow described above, it is possible to obtain a self-contained, high-per{\-}form{\-}ance solver that can be easily deployed on embedded hardware starting from a description of the OCP in a high-level language.

We remark that model equations and other nonlinear functions are called from \texttt{acados} in a completely language-agnostic way: \software{acados} is at no point aware of which modeling tool is being used.
One benefit is that this facilitates self-written models (in \texttt{C/C++}), which are also completely compatible with \texttt{acados}.
However, this is more involved, since in the case of \software{CasADi} functions, memory allocation and matrix format conversions are taken into account automatically by the \software{CasADi} functions wrapper in \software{acados}.

\section{Numerical Results} \label{sec:numerical_results}

This section consists of a few numerical experiments with~\software{acados} and comparisons to other embedded optimization software packages.
We discuss performance on the nonlinear chain-of-masses problem, we present one open-loop example with different Hessian approximations, and show one closed-loop engine control experiment on an embedded platform.

\subsection{Case Study 1: Chain of Masses} \label{sec:chain}

As a benchmarking problem, we take the chain-of-masses problem as presented in~\cite{Wirsching2006}. The control objective is to stabilize a chain of masses with nonlinear interaction between them. For a full description of the system, we refer to the appendix.
The system is useful as a benchmark in the sense that the problem is simple enough to understand intuitively, yet complicated enough to get non-trivial results from a range of different solvers.
Also, by increasing the number of masses, one could compare behavior for different numbers of states easily, without changing much code.

\subsubsection*{Closed-loop experiments}

In closed-loop, an MPC controller repeatedly (approximately) solves OCP~\eqref{eq:chain_of_masses_ocp}.
The first control $u_0$ is passed to the dynamic system under control and a new initial state $\overline{x}_0$ is obtained.
Here, we simulate the system by using a more accurate integrator than the one in OCP~\eqref{eq:chain_of_masses_ocp}, namely the Dormand-Prince method, as implemented in the \Matlab{} routine \software{ode45}.

We introduce one disturbance into the closed-loop system, similar as in~\cite{Wirsching2006}: in the beginning of the simulation, we start from a horizontal configuration of the chain of masses.
Around the midpoint of the simulation, we override the closed loop control with a constant $u_\mathrm{d} = [-1, 1, 1]^\top$.
After one second of simulation time, the controller takes over again.

We compare the following solvers with each other for this particular closed-loop setup:

\begin{itemize}
	\item \software{IPOPT}~\cite{Waechter2006}.
	As a solver \emph{not} targeting embedded devices specifically, we use it as a baseline to compare against.
	\item \software{FalcOPT}~\cite{Torrisi2018}.
	A projected gradient descent method tailored for NMPC.
	\item \software{VIATOC}~\cite{Kalmari2015}.
	A gradient projection method for MPC that only allows linear inequality constraints.
	\item \software{ACADO Code Generation} tool~\cite{Houska2011}.
	Generates SQP-based solvers.
	\item \software{GRAMPC}~\cite{Englert2019}.
	An embedded Augmented Lagrangian-based solver.
	\item \software{acados}.
	Framework presented in the current paper.
\end{itemize}
The tuning parameters for the different solvers are listed in Table~\ref{tab:tuning_parameters}.

\begin{table}[h]
	\centering
	\caption{\textsc{Solver options for the different solvers in Case Study 1.}}
	\label{tab:tuning_parameters}
	\begin{tabular}{|c|l|}
		\hline Solver & Solver options \\ \hline
		\software{IPOPT} & Called through \software{CasADi}, default parameters \\ 
		\software{FalcOPT} & Tolerance (\texttt{eps}): 0.1, maximum number of iterations (\texttt{maxIt}): 100 \\ 
		\software{VIATOC} & Maximum number of iterations: 20 \\ 
		\software{ACADO} & RTI solver, Full condensing, QP solver \software{qpOASES} \\ 
		\software{GRAMPC} & Parameters chosen as in~\cite{Englert2019}: max. number of augm. Lagrange iterations: 5 \\
		\software{acados} & \texttt{SQP\_RTI} solver, QP solver \software{HPIPM}, partial condensing horizon of 5 \\ \hline 
	\end{tabular} 
\end{table}

In order to compare the quality of the closed-loop solutions, we use the notion of distance-from-reference (DR), which is an approximation of the integrated cost along closed-loop trajectories:
\begin{equation*}
\mathrm{DR}_{(\cdot), n} = \sum_{i=0}^{n} \begin{bmatrix}x_{(\cdot),i}-x_\mathrm{ref} \\ u_{(\cdot),i}-u_\mathrm{ref}\end{bmatrix}^\top \begin{bmatrix} Q & \zeros \\ \zeros & R
\end{bmatrix} \begin{bmatrix}x_{(\cdot),i}-x_\mathrm{ref} \\ u_{(\cdot),i}-u_\mathrm{ref}\end{bmatrix}.
\end{equation*}

To compare the different solvers, we plot the relative cumulative sub-optimality (RCSO), relative to a fully converged solution, in this case, the \software{IPOPT} solution, which reads as
\begin{equation*}
\mathrm{RCSO}_{(\cdot),n} = \left|\frac{\mathrm{DR}_{(\cdot),n} - \mathrm{DR}_{\software{ipopt},n}}{\mathrm{DR}_{\software{ipopt},n}}\right|,
\end{equation*}
where $n = 0,\ldots, 300$ denotes the time step in our simulation.
We show a comparison in  Table~\ref{tab:nlchainclosedloop}.
The results for \software{ACADO} and \software{acados} are \emph{exactly} the same, as they implement the same real-time algorithm, with both being very close to the reference solution from \software{IPOPT}.
The solvers \software{GRAMPC}, \software{VIATOC} and \software{FalcOPT}, being based on first-order methods, are further away from the \software{IPOPT} solution.
These findings are consistent with previously published work by other authors, see \cite{Englert2019}.

\begin{table}
	\centering
	\caption{Relative suboptimality at the end of the simulation of the hanging chain with $M=5$ and $N=40$.
	First-order methods \software{VIATOC}, \software{GRAMPC} and \software{FalcOPT} were tuned to perform similarly.
	The algorithms chosen in \software{ACADO} and \software{acados} are identical, hence the results are identical.}
	\label{tab:nlchainclosedloop}
	\begin{tabular}{|c|c|}
		\hline Solver name & RCSO \\
		\hline \texttt{IPOPT} & 0.00e+00 \\
		\hline \texttt{FalcOPT} & 3.170e-01 \\
		\hline \texttt{VIATOC} & 4.74e-03 \\
		\hline \texttt{ACADO} & 1.01e-04 \\
		\hline \texttt{GRAMPC} & 7.17e-02 \\ 
		\hline \texttt{acados} & 1.01e-04 \\ \hline
	\end{tabular}
\end{table}

\begin{figure}
	\centering
	\includegraphics[width=\linewidth]{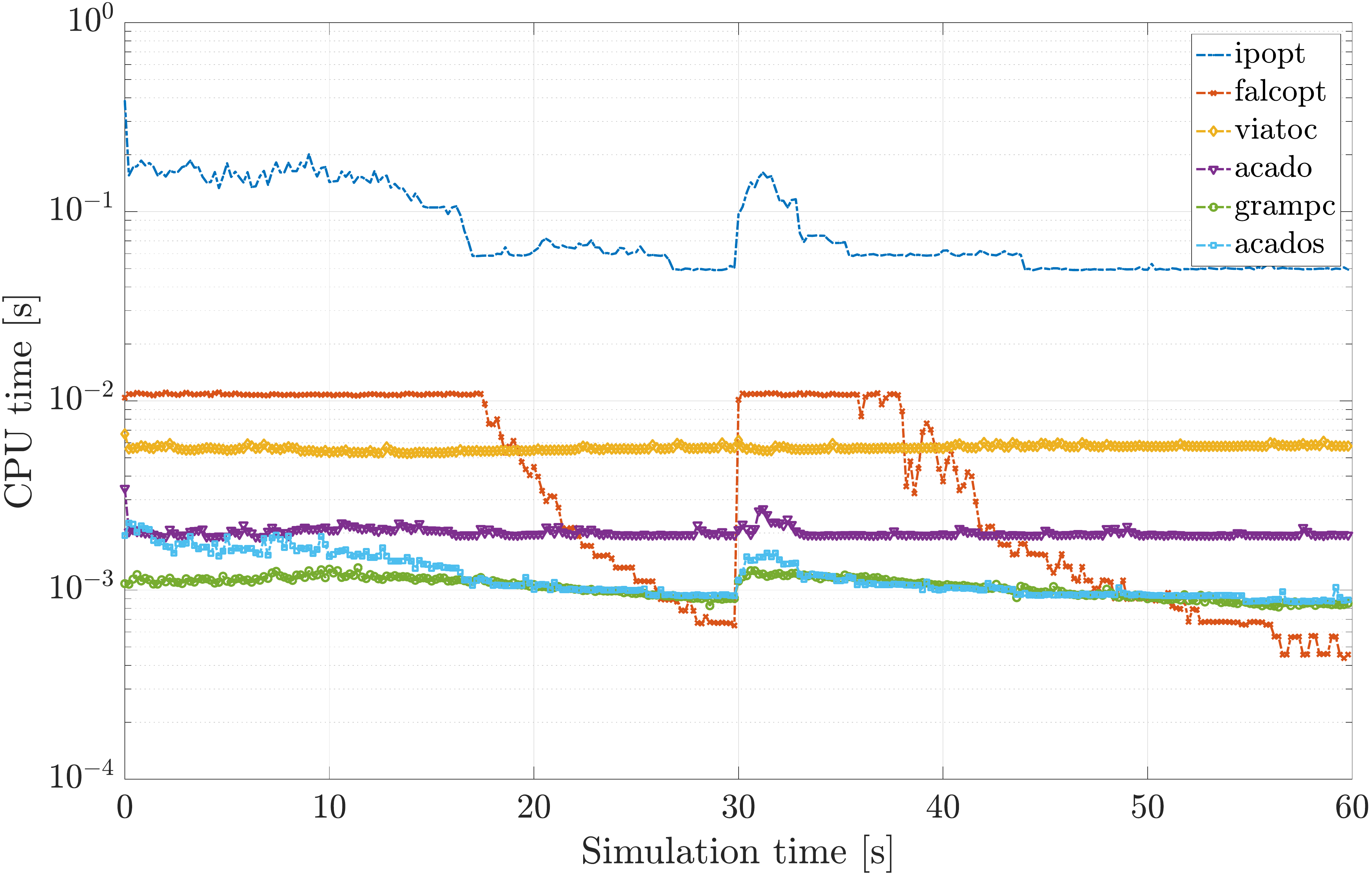}
	\caption{Computational time for each iteration of the closed loop simulation, averaged over 10 runs.}
	\label{fig:nlchaintimings}
\end{figure}

We have a look at the computational performance along the closed-loop trajectories in Figure~\ref{fig:nlchaintimings}.
\software{GRAMPC}, \software{ACADO}, \software{VIATOC} and \software{acados} produce consistent timings throughout the entire experiment, even when the disturbance occurs.
This is a beneficial property for embedded solvers, as they often have a fixed time deadline, being part of a larger control application.
\software{GRAMPC} and \software{acados} produce solutions at almost the same speed, both approximately a factor 2 faster than \software{ACADO} which is in turn a factor 2-3 faster than \software{VIATOC}.
Near the equilibrium, \software{FalcOPT} takes the shortest computation time, as it needs to perform only a few gradient steps per iteration.
\software{IPOPT} is included as a baseline for comparison to non-embedded solvers.
The timings are summarized in Table~\ref{tab:nlchaintimings}.

\begin{table}
	\centering
	\caption{\textsc{Computation times for the closed-loop experiments on a chain of masses (cf.
	Figure~\ref{fig:nlchaintimings})}}
	\begin{tabular}{|c|c|c|c|}
		\hline
		comp. time per iteration $(\mathrm{ms})$ & median & minimum & maximum \\ 
		\hline \hline
		\software{IPOPT} & 59.84 & 49.06 & 384.90 \\ 
		\hline
		\software{FalcOPT} & 4.36 & 0.44 & 11.10 \\
		\hline
		\software{VIATOC} & 5.63 & 5.27 & 6.68 \\
		\hline
		\software{ACADO} & 1.97 & 1.90 & 3.45 \\
		\hline
		\software{GRAMPC} & 1.06 & 0.81 & 1.31 \\
		\hline
		\software{acados} & 1.05 & 0.87 & 2.23 \\
		\hline 
	\end{tabular}
	\label{tab:nlchaintimings}
\end{table}

\begin{figure}
	\centering
	\includegraphics[width=0.7\linewidth]{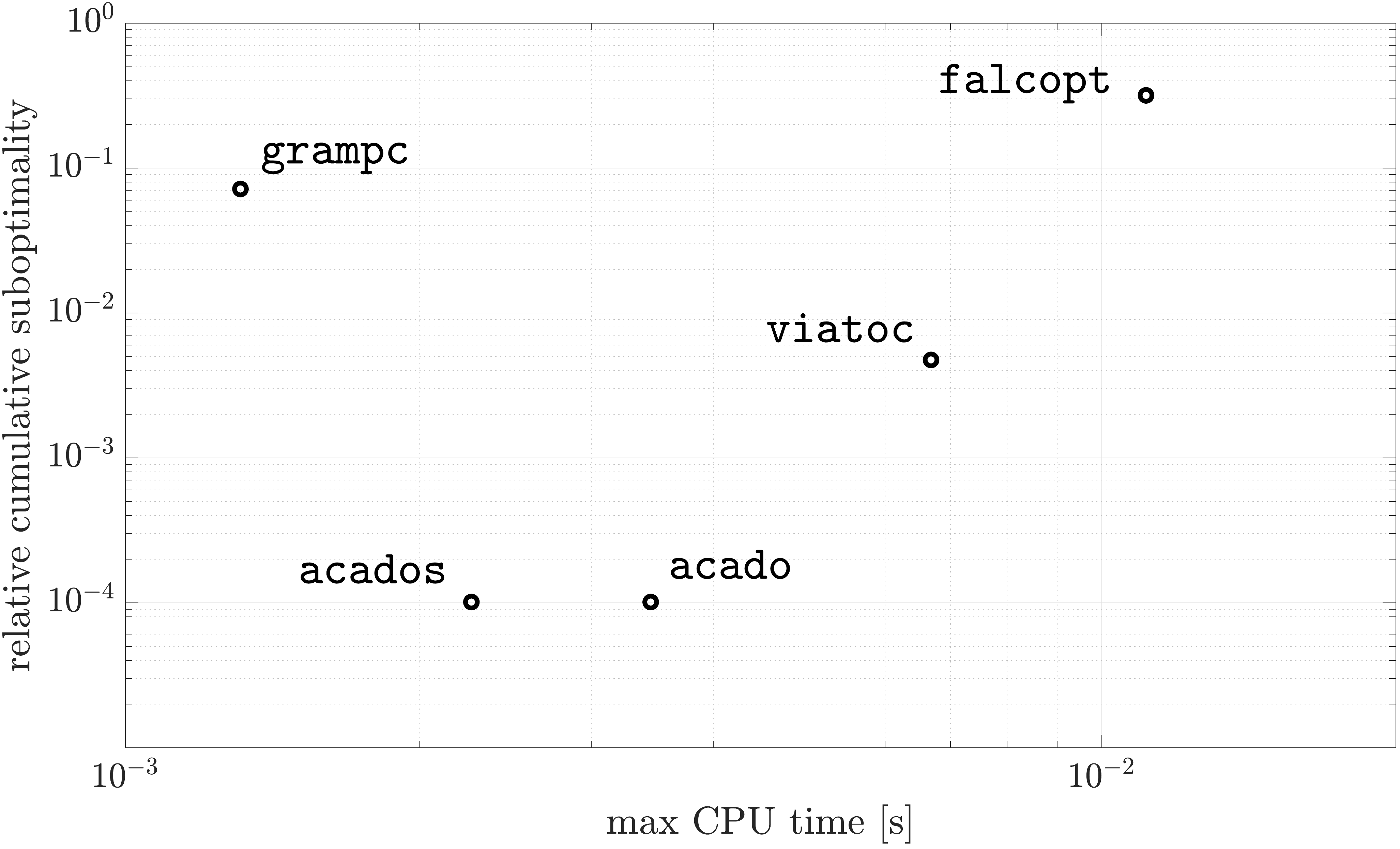}
	\caption{Trade-off between sub-optimality (Table~\ref{tab:nlchainclosedloop}) and computation time (Figure~\ref{fig:nlchaintimings}).
	We see that \software{acados} and \software{GRAMPC} lie on the Pareto-optimal front.}
	\label{fig:nlchainparetofront}
\end{figure}

Of course, any solver can trade off sub-optimality for computation time.
To get the full picture, we plot both measures against each other in Figure~\ref{fig:nlchainparetofront}: we look at relative cumulative sub-optimality over the entire length of the experiment, versus worst-case computation times.
By this comparison, we see that \software{acados} and \software{GRAMPC} are on the Pareto-optimal front: although \software{acados} is a factor 1000 less suboptimal than \software{GRAMPC}, the computational cost is higher.
By the median computation times, \software{acados} is faster (see Table~\ref{tab:nlchaintimings}).

\subsection{Case study 2: Hessian regularization}

In Section~\ref{sec:embedded_sqp}, we briefly mentioned the impact of Hessian regularization on SQP methods.
In this case study, we compare the convergence of exact-Hessian based SQP with three different Hessian regularizations, on a simple control problem, namely a cart-pole swingup, which is described in the appendix.
Note that in this case study, as opposed to the others, we do not perform closed-loop experiments but solve one optimal control problem up until convergence.

\subsubsection{Exact-Hessian based SQP}

We solve OCP~\eqref{eq:OCP_pendulum} with SQP, where we use the exact Hessian of the Lagrangian.
In the notation of~\eqref{eq:OCP_QP_subproblem}:

\begin{align*}
H_k &= \begin{bmatrix}
Q & \zeros \\ \zeros & R
\end{bmatrix} + \sum_{i=0}^{n_x} \pi_{k,i} \nabla^2_{(x,u)} \phi^x_i(x_k, u_k), \quad k=0,\ldots,N-1 \\
H_N &= Q,
\end{align*}
where $\pi_{k,i}$ are the Lagrange multipliers associated with the dynamic equality constraints.

In some cases, the non-convexity of the dynamic equations gives rise to an indefinite Hessian matrix.
We apply the $\mathrm{project}(\cdot)$ and $\mathrm{mirror}(\cdot)$ regularizations, as well as the convexification method previously mentioned in Section~\ref{sec:convex_hessians}.
All are implemented as modules in the \software{acados} framework.

\begin{figure}
	\centering
	\includegraphics[width=0.7\linewidth]{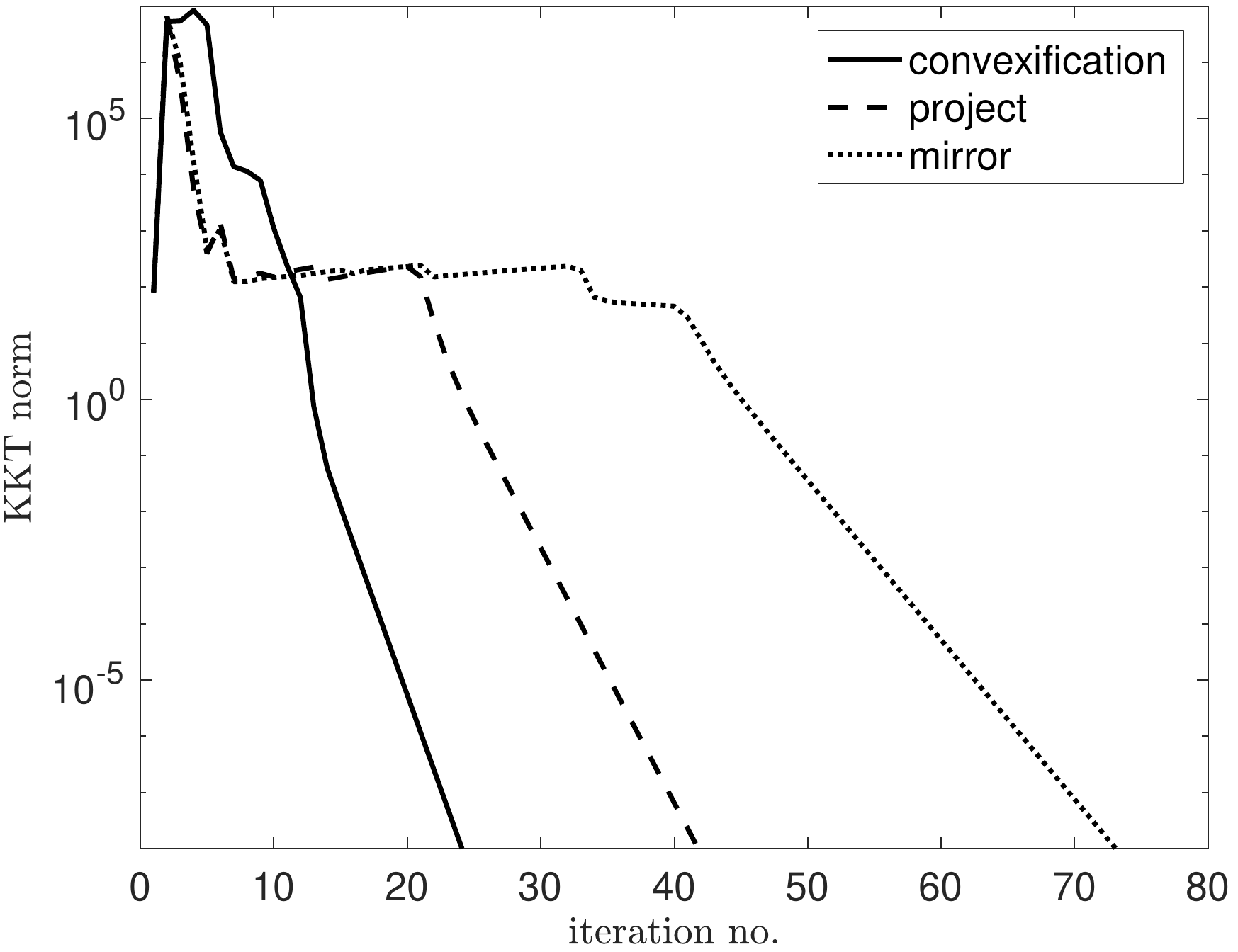}
	\caption{Convergence comparison of exact-Hessian based SQP with three different regularization strategies.}
	\label{fig:regularization_convergence}
\end{figure}

We compare the convergence of the SQP iterates obtained using the three different regularization methods.
For each SQP variant, we start the SQP iterations from an initialization point with zeros for all states except a linearly decreasing initialization for the angle, from $\pi$ to zero.
The result can be seen  in Figure~\ref{fig:regularization_convergence}.
The structure-exploiting convexification converges almost twice as fast as the projection regularization, and is in turn much faster than the mirroring regularization.
Intuitively, this makes sense, as mirroring is `blocking' directions associated with large negative eigenvalues, by introducing large positive eigenvalues in those directions.
This prevents the solver from taking larger steps\footnote{This is also the approach followed by the algorithms obtained with the \software{ACADO Code Generation} tool.}.
In turn, the structure-exploiting regularization is faster than merely projecting the eigenvalues on the positive definite cone, because it is redistributing convexity among all stages, and thus needs less regularization overall.

\begin{table}
	\centering
	\caption{\textsc{Exact-Hessian based SQP: computation times}}
	\label{tab:regularization_timings}
	\begin{tabular}{|c||c|c|c|}
		\hline regularization & convexification & project & mirror \\
		\hline avg iteration (ms) & 2.540 & 2.303 & 2.264 \\
		\hline total time (ms) & 66.034 & 103.65 & 174.32 \\ \hline
	\end{tabular}
\end{table}

It must be said that the convexification method is quite a bit more involved than the other two regularization schemes.
However, by using the optimized linear algebra routines of \software{BLASFEO}, we implemented the convexification method such that it is only slightly more expensive per iteration than the basic regularization methods, see Table~\ref{tab:regularization_timings}, but much less computationally expensive overall.
Thus, the Hessian convexification method allows us to perform exact-Hessian based NMPC online, with better performance than state-of-the-art methods.

\subsection{Case study 3: Hardware-in-the-loop experiments for an engine control application}

As a last case study, we discuss the performance of \software{acados} on an embedded platform, namely the dSPACE MicroAutoboxII~\cite{dspace2006}.
It is an industrial computing platform that is used in the car industry.
It features a 900 MHz PowerPC processor (IBM PPC 750GL) with 16MB of main memory.
The control application that we focus on is engine control, with the engine model as presented in~\cite{Albin2017}, which we will briefly reproduce in the appendix.


The control objective is to track a boost pressure signal, where the boost pressure is given by $y_p(x) := \Pi_{c,\mathrm{lp}} \cdot \Pi_{c,\mathrm{hp}}$ (see appendix).
To this end, we solve an OCP arising from a multiple shooting formulation with the Gauss-Legendre method of order 6 with sampling time $0.05 \, \mathrm{s}$ and $N=20$ shooting intervals.
The DAE simulation functions are denoted by $\phi$.
Let $r(x, u) = [y_p(x); x; u]$ and $r_N(x) = [y_p(x); x]$.
The OCP then reads as
\begin{mini}
	{\substack{x_0, \ldots, x_N, \\ z_0, \ldots, z_{N-1} \\ u_0, \ldots, u_{N-1}}}
	{\sum_{k=0}^{N-1} \|r(x_k, u_k) - y_{r,k}\|_W^2 + \|r_N(x_N) - y_{r,N}\|_{W_N}^2}
	{\label{eq:engine_control}}
	{}
	\addConstraint{x_0}{= \overline{x}_0}{}
	\addConstraint{\begin{bmatrix}x_{k+1} \\ z_k \end{bmatrix}}{= \phi(x_k, u_k), \quad}{k=0,\ldots,N-1}
	\addConstraint{0}{\leq u_k \leq 100,}{k=0,\ldots,N-1}
	\addConstraint{0.5}{\leq \Pi_{c,\mathrm{lp}, k} \leq 1.757,\quad }{k=1,\ldots,N}
	\addConstraint{0.5}{\leq \Pi_{c,\mathrm{hp}, k} \leq 2.125,}{k=1,\ldots,N.}
\end{mini}

Constraints on $\Pi_{c,\mathrm{lp}}, \Pi_{t,\mathrm{hp}}$ are included to prevent damage to the compressor. The exact values of the weight matrices and the reference vectors can be found in the appendix.

\begin{figure}
	\centering
	\includegraphics[width=\linewidth]{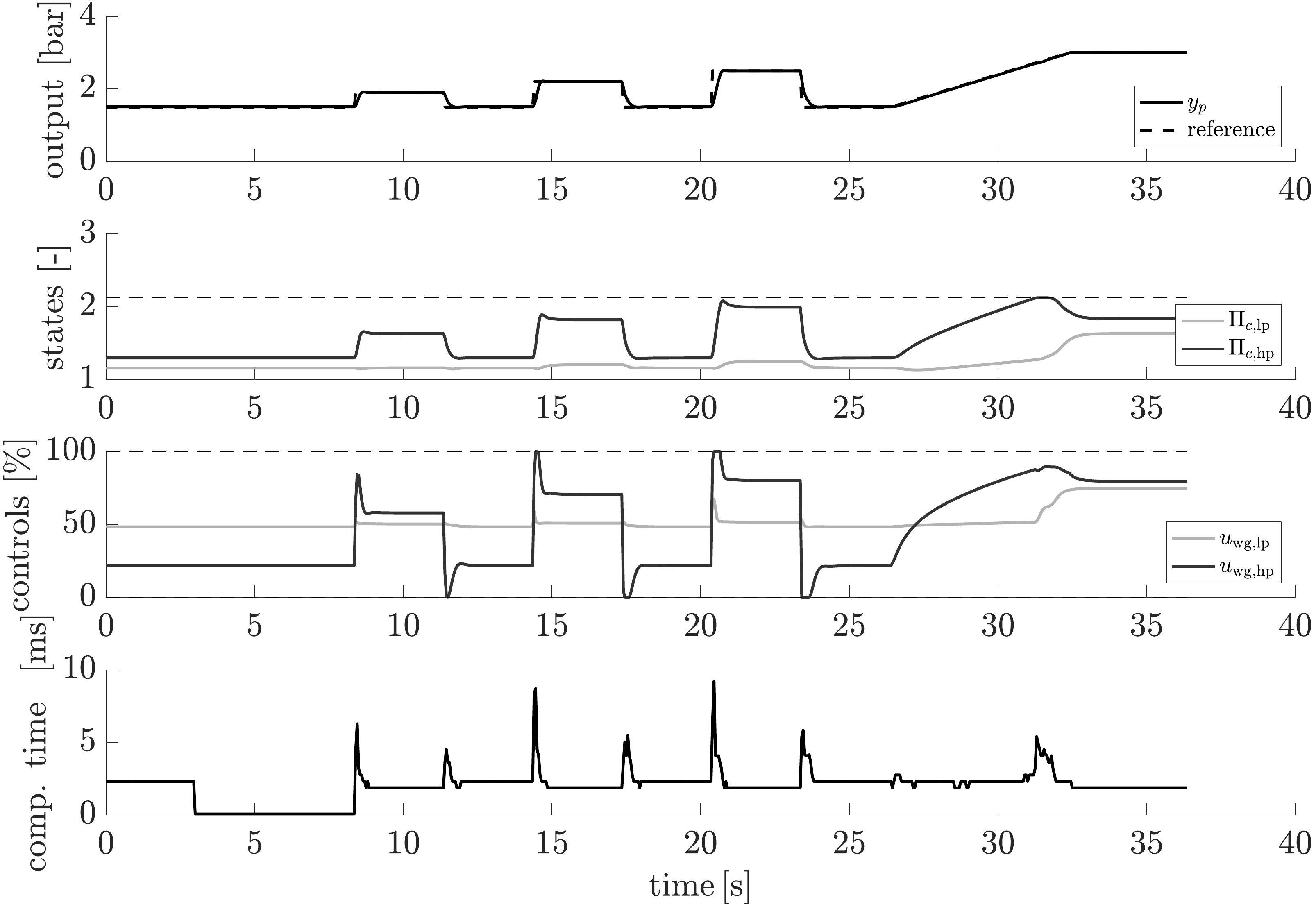}
	\caption{Closed-loop simulation of the engine control task with steps in the reference boost pressure.
		Simulations are carried out on the dSPACE MicroAutoboxII platform at a clock speed of $900\,\mathrm{MHz}$.}
	\label{fig:enginetracking}
\end{figure}

We repeatedly solve OCP~\eqref{eq:engine_control} approximately by performing real time iterations.
As an underlying QP solver, we use~\software{HPIPM}.
When ran in closed loop on the dSPACE MicroAutoboxII, the results can be seen in Figure~\ref{fig:enginetracking}.
Control bounds and state bounds become active at some point in the simulation, for the high-pressure stage.
The reference is tracked closely and without oscillations, which have been observed when linear-quadratic MPC is used~\cite{Albin2017}.
As for the computation times, it is interesting to note that there are spikes everywhere where a jump occurs or a constraint becomes (in)active.
The computation times close to the solution (i.e.
at the beginning of the simulation) drop to almost zero.
Overall, the maximum computation time remains under $10\,\mathrm{ms}$, which is 5x faster than the sampling time of the system ($50\,\mathrm{ms}$).
We remark that the computation times obtained with the dSPACE MicroAutoboxII, for this HIL experiment, are about three times slower than a desktop computer with a 2.5GHz Intel Core i7-4870HQ processor.

\section{Conclusion and outlook}
\label{sec:conclusion}

In this article, we presented \software{acados}, a new software package for embedded optimization.
It is free and open-source software that facilitates rapid testing and deployment of (N)MPC algorithms on embedded hardware platforms.
For ease of use, we offer interfaces with higher-level languages such as \textsc{Matlab} and \software{Python}.

Among many features that state-of-the-art NMPC algorithms require, a couple of new features that are not present in any other software package is the convexification procedure of Section~\ref{sec:convex_hessians}, allowing the use of exact-Hessian based SQP methods in real-time, and the SCQP Hessian approximation.
Additionally, the structure exploiting GNSF-IRK integrator has the potential to speed up the simulation and sensitivity propagation tasks within an NMPC scheme.
Furthermore, \software{acados} features partial condensing, different state and control dimensions per multiple shooting stage, the use of \software{BLASFEO} as a linear algebra backend and facilities for using \software{CasADi} as a modeling language.

The software is shown to be embeddable, by numerical experiments on the dSPACE MicroAutoboxII industrial computer, resulting in computation times in the millisecond range for a non-trivial NMPC problem.
Furthermore, it is shown to be fast, by comparison to other embedded optimization packages.

\software{acados} is an ongoing endeavor.
Future work includes extending interoperability with \software{Simulink} for easier deployment on embedded systems, and adding features for nonlinear interior-point methods, as well as other SQP-based methods like multi-level iterations.


\bibliographystyle{spmpsci}        
\bibliography{syscop.bib}           

\begin{thebibliography}{10}
\providecommand{\url}[1]{{#1}}
\providecommand{\urlprefix}{URL }
\expandafter\ifx\csname urlstyle\endcsname\relax
  \providecommand{\doi}[1]{DOI~\discretionary{}{}{}#1}\else
  \providecommand{\doi}{DOI~\discretionary{}{}{}\begingroup
  \urlstyle{rm}\Url}\fi

\bibitem{Albersmeyer2010b}
Albersmeyer, J.: Adjoint-based algorithms and numerical methods for sensitivity
  generation and optimization of large scale dynamic systems.
\newblock Ph.D. thesis, University of Heidelberg
\newblock  (2010)

\bibitem{Albersmeyer2010}
Albersmeyer, J., Diehl, M.: The lifted {N}ewton method and its application in
  optimization.
\newblock SIAM Journal on Optimization 20(3), 1655--1684
\newblock  (2010)

\bibitem{Albin2017}
Albin, T., Ritter, D., Liberda, N., Quirynen, R., Diehl, M.: In-vehicle
  realization of nonlinear {MPC} for gasoline two-stage turbocharging airpath
  control.
\newblock IEEE Transactions on Control Systems Technology pp. 1--13
\newblock  (2017)

\bibitem{Andersson2018}
Andersson, J.A.E., Gillis, J., Horn, G., Rawlings, J.B., Diehl, M.: {CasADi}: a
  software framework for nonlinear optimization and optimal control.
\newblock Mathematical Programming Computation
\newblock  (2018)

\bibitem{Andersson2018a}
Andersson, J.A.E., Rawlings, J.B.: Sensitivity analysis for nonlinear
  programming in casadi.
\newblock In: Proceedings of the IFAC Conference on Nonlinear Model Predictive
  Control (NMPC)
\newblock  (2018)

\bibitem{Axehill2015}
Axehill, D.: Controlling the level of sparsity in {MPC}.
\newblock Systems \& Control Letters 76, 1--7
\newblock  (2015)

\bibitem{Bemporad2005}
Bemporad, A.: {H}ybrid {T}oolbox for Matlab
\newblock  (2003)

\bibitem{Bemporad1999}
Bemporad, A., Borrelli, F., Morari, M.: The explicit solution of constrained
  {LP}-based receding horizon control.
\newblock In: Proceedings of the IEEE Conference on Decision and Control (CDC).
  Sydney, Australia
\newblock  (1999)

\bibitem{Bock1987}
Bock, H.: {R}andwertproblemmethoden zur {P}arameteridentifizierung in
  {S}ystemen nichtlinearer {D}ifferentialgleichungen, \emph{Bonner
  Ma\-the\-ma\-tische Schriften}, vol. 183.
\newblock Universit\"at Bonn, Bonn
\newblock  (1987)

\bibitem{Bock1983}
Bock, H.G.: Recent advances in parameter identification techniques for {ODE}.
\newblock In: Numerical Treatment of Inverse Problems in Differential and
  Integral Equations, pp. 95--121. Birk\-h\"au\-ser
\newblock  (1983)

\bibitem{Bock2005}
Bock, H.G., Diehl, M., K\"uhl, P., Kostina, E., Schl\"oder, J.P., Wirsching,
  L.: Numerical methods for efficient and fast nonlinear model predictive
  control.
\newblock In: Proceedings of ''Int. Workshop on assessment and future
  directions of Nonlinear Model Predictive Control''. Springer
\newblock  (2005)

\bibitem{Bock1984}
Bock, H.G., Plitt, K.J.: A multiple shooting algorithm for direct solution of
  optimal control problems.
\newblock In: Proceedings of the IFAC World Congress, pp. 242--247. Pergamon
  Press
\newblock  (1984)

\bibitem{Cheshmi2020}
Cheshmi, K., Kaufman, D.M., Kamil, S., Dehnavi, M.M.: {NASOQ}: Numerically
  accurate sparsity-oriented qp solver.
\newblock
\newblock In: ACM Transactions on Graphics, vol.~39

\bibitem{Chiang2016a}
Chiang, N.Y., Hang, R., Zavala, V.M.: An augmented lagrangian filter method for
  real-time embedded optimization.
\newblock IEEE Transactions on Automatic Control 62(12), 6110--6121
\newblock  (2017)

\bibitem{Cimini2017}
Cimini, G., Bemporad, A.: Exact complexity certification of active-set methods
  for quadratic programming.
\newblock IEEE Transactions on Automatic Control 62(12), 6094--6109
\newblock  (2017)

\bibitem{Deng2018}
Deng, H., Ohtsuka, T.: A highly parallelizable newton-type method for nonlinear
  model predictive control.
\newblock In: Proceedings of the IFAC Conference on Nonlinear Model Predictive
  Control (NMPC)
\newblock  (2018)

\bibitem{Cairano2013}
{Di Cairano}, S., Brand, M., Bortoff, S.A.: Projection-free parallel quadratic
  programming for linear model predictive control.
\newblock International Journal of Control 86(8), 1367--1385
\newblock  (2013)

\bibitem{Diehl2002b}
Diehl, M., Bock, H.G., Schl\"oder, J.P., Findeisen, R., Nagy, Z., Allg\"ower,
  F.: Real-time optimization and nonlinear model predictive control of
  processes governed by differential-algebraic equations.
\newblock Journal of Process Control 12(4), 577--585
\newblock  (2002)

\bibitem{Diehl2009c}
Diehl, M., Ferreau, H.J., Haverbeke, N.: Efficient numerical methods for
  nonlinear {MPC} and moving horizon estimation.
\newblock In: L.~Magni, M.~Raimondo, F.~Allg\"ower (eds.) Nonlinear model
  predictive control, \emph{Lecture Notes in Control and Information Sciences},
  vol. 384, pp. 391--417. Springer
\newblock  (2009)

\bibitem{Domahidi2013a}
Domahidi, A., Chu, E., Boyd, S.: {ECOS}: An {SOCP} solver for embedded systems.
\newblock In: Proceedings of the European Control Conference (ECC), pp.
  3071--3076. IEEE
\newblock  (2013)

\bibitem{Domahidi2012}
Domahidi, A., Zgraggen, A., Zeilinger, M.N., Morari, M., Jones, C.N.: Efficient
  interior point methods for multistage problems arising in receding horizon
  control.
\newblock In: Proceedings of the IEEE Conference on Decision and Control (CDC),
  pp. 668--674. Maui, HI, USA
\newblock  (2012)

\bibitem{dspace2006}
d{SPACE}: {H}omepage.
\newblock http://www.dspace.com
\newblock  (2006)

\bibitem{Englert2019}
Englert, T., Völz, A., Mesmer, F., Rhein, S., Graichen, K.: A software
  framework for embedded nonlinear model predictive control using a
  gradient-based augmented lagrangian approach ({GRAMPC}).
\newblock Optimization and Engineering 20(3), 769--809
\newblock  (2019)

\bibitem{Ferreau2017}
Ferreau, H.J., Almer, S., Verschueren, R., Diehl, M., Frick, D., Domahidi, A.,
  Jerez, J.L., Stathopoulos, G., Jones, C.: Embedded optimization methods for
  industrial automatic control.
\newblock In: Proceedings of the IFAC World Congress
\newblock  (2017)

\bibitem{Ferreau2014}
Ferreau, H.J., Kirches, C., Potschka, A., Bock, H.G., Diehl, M.: {qpOASES}: a
  parametric active-set algorithm for quadratic programming.
\newblock Mathematical Programming Computation 6(4), 327--363
\newblock  (2014)

\bibitem{Franke1997}
Franke, R., Arnold, E.: Computer Intensive Methods in Control and Signal
  Processing, chap. Applying new numerical algorithms to the solution of
  discrete-time optimal control problems, pp. 105--117.
\newblock Springer
\newblock  (1997)

\bibitem{Frasch2015}
Frasch, J.V., Sager, S., Diehl, M.: A parallel quadratic programming method for
  dynamic optimization problems.
\newblock Mathematical Programming Computations 7(3), 289--329
\newblock  (2015)

\bibitem{Frey2020}
Frey, J., Cairano, S.D., Quirynen, R.: {A}ctive-{S}et based {I}nexact
  {I}nterior {P}oint {QP} {S}olver for {M}odel {P}redictive {C}ontrol.
\newblock In: Proceedings of the IFAC World Congress
\newblock  (2020)

\bibitem{Frey2019}
Frey, J., Quirynen, R., Kouzoupis, D., Frison, G., Geisler, J., Schild, A.,
  Diehl, M.: Detecting and exploiting {G}eneralized {N}onlinear {S}tatic
  {F}eedback structures in {DAE} systems for {MPC}.
\newblock In: Proceedings of the European Control Conference (ECC)
\newblock  (2019)

\bibitem{Frison2015a}
Frison, G.: Algorithms and methods for high-performance model predictive
  control.
\newblock Ph.D. thesis, Technical University of Denmark (DTU)
\newblock  (2015)

\bibitem{Frison2020a}
Frison, G., Diehl, M.: {HPIPM}: a high-performance quadratic programming
  framework for model predictive control.
\newblock In: Proceedings of the IFAC World Congress. Berlin, Germany
\newblock  (2020)

\bibitem{Frison2018}
Frison, G., Kouzoupis, D., Sartor, T., Zanelli, A., Diehl, M.: {BLASFEO}: Basic
  linear algebra subroutines for embedded optimization.
\newblock ACM Transactions on Mathematical Software (TOMS) 44(4), 42:1--42:30
\newblock  (2018)

\bibitem{Frison2017}
Frison, G., Quirynen, R., Zanelli, A., Diehl, M., J{\o}rgensen, J.B.: Hardware
  tailored linear algebra for implicit integrators in embedded {NMPC}.
\newblock In: Proceedings of the IFAC World Congress
\newblock  (2017)

\bibitem{Frison2014}
Frison, G., Sorensen, H.B., Dammann, B., J{\o}rgensen, J.B.: High-performance
  small-scale solvers for linear model predictive control.
\newblock In: Proceedings of the European Control Conference (ECC), pp.
  128--133
\newblock  (2014)

\bibitem{Gertz2003}
Gertz, E.M., Wright, S.J.: Object-oriented software for quadratic programming.
\newblock ACM Transactions on Mathematical Software 29(1), 58--81
\newblock  (2003)

\bibitem{Giftthaler2018}
Giftthaler, M., Neunert, M., {St\"auble}, M., Buchli, J.: The {Control Toolbox}
  - an open-source {C++} library for robotics, optimal and model predictive
  control.
\newblock In: IEEE International Conference on Simulation, Modeling, and
  Programming for Autonomous Robots (SIMPAR)
\newblock  (2018)

\bibitem{Griewank2000}
Griewank, A.: {E}valuating {D}erivatives, {P}rinciples and {T}echniques of
  {A}lgorithmic {D}ifferentiation.
\newblock No.~19 in Frontiers in Appl. Math. {SIAM}, Philadelphia
\newblock  (2000)

\bibitem{Gros2006}
Gros, S., Srinivasan, B., Bonvin, D.: {R}obust predictive control based on
  neighboring extremals.
\newblock Journal of Process Control 16, 243--253
\newblock  (2006)

\bibitem{Gros2016}
Gros, S., Zanon, M., Quirynen, R., Bemporad, A., Diehl, M.: From linear to
  nonlinear {MPC}: bridging the gap via the real-time iteration.
\newblock International Journal of Control
\newblock  (2016)

\bibitem{Gruene2017}
Gr{\"u}ne, L., Pannek, J.: Nonlinear Model Predictive Control, second edition
  edn.
\newblock Springer-Verlag
\newblock  (2017)

\bibitem{Hairer1993}
Hairer, E., N{\o}rsett, S., Wanner, G.: {S}olving {O}rdinary {D}ifferential
  {E}quations {I}, 2nd edn.
\newblock Springer Series in Computational Mathematics. Springer, Berlin
\newblock  (1993)

\bibitem{Hairer1991}
Hairer, E., Wanner, G.: {S}olving {O}rdinary {D}ifferential {E}quations {II} --
  {S}tiff and {D}ifferential-{A}lgebraic {P}roblems, 2nd edn.
\newblock Springer, Berlin Heidelberg
\newblock  (1991)

\bibitem{Hehn2011}
Hehn, M., D'Andrea, R.: A flying inverted pendulum.
\newblock In: IEEE International Conference on Robotics and Automation,
\newblock pp. 763--770

\bibitem{Herceg2013}
Herceg, M., Kvasnica, M., Jones, C., Morari, M.: {Multi-Parametric Toolbox
  3.0}.
\newblock In: Proc.~of the European Control Conference, pp. 502--510. Z\"urich,
  Switzerland
\newblock  (2013).
\newblock \url{http://control.ee.ethz.ch/~mpt}

\bibitem{Hermans2020}
Hermans, B., Themelis, A., Patrinos, P.:
\newblock Qpalm: A proximal augmented lagrangian method for nonconvex quadratic
  programs

\bibitem{Houska2011}
Houska, B., Ferreau, H.J., Diehl, M.: An auto-generated real-time iteration
  algorithm for nonlinear {MPC} in the microsecond range.
\newblock Automatica 47(10), 2279--2285
\newblock  (2011)

\bibitem{Kalmari2015}
Kalmari, J., Backman, J., Visala, A.: A toolkit for nonlinear model predictive
  control using gradient projection and code generation.
\newblock Control Engineering Practice 39, 56--66
\newblock  (2015)

\bibitem{Kapernick2014}
K{\"a}pernick, B., Graichen, K.: The gradient based nonlinear model predicitive
  control software {GRAMPC}.
\newblock In: Proceedings of the European Control Conference (ECC)
\newblock  (2014)

\bibitem{Katliar2020}
Katliar, M.: Optimal control of motion simulators.
\newblock Ph.D. thesis, Albert-Ludwigs-Universit{\"a}t Freiburg
\newblock  (2020)

\bibitem{Khusainov2017}
Khusainov, B., Kerrigan, E.C., Suardi, A., Constantinides, G.A.: Nonlinear
  predictive control on a heterogeneous computing platform.
\newblock In: Proceedings of the IFAC World Congress
\newblock  (2017)

\bibitem{Kirches2012b}
Kirches, C., Bock, H.G., Schl\"oder, J.P., Sager, S.: Complementary condensing
  for the direct multiple shooting method.
\newblock In: Modeling, Simulation and Optimization of Complex Processes, pp.
  195--206. Springer Berlin Heidelberg
\newblock  (2012)

\bibitem{Kouzoupis2015}
Kouzoupis, D., Ferreau, H.J., Peyrl, H., Diehl, M.: First-order methods in
  embedded nonlinear model predictive control.
\newblock In: Proceedings of the European Control Conference (ECC), pp.
  2617--2622
\newblock  (2015)

\bibitem{Kouzoupis2015b}
Kouzoupis, D., Quirynen, R., Frasch, J.V., Diehl, M.: Block condensing for fast
  nonlinear {MPC} with the dual {N}ewton strategy.
\newblock In: Proceedings of the IFAC Conference on Nonlinear Model Predictive
  Control (NMPC), vol.~48, pp. 26--31
\newblock  (2015)

\bibitem{Kvamme2014}
Kvamme, S.: {DuQuad} {W}ebpage.
\newblock \url{http://sverrkva.github.io/duquad/}
\newblock  (2014)

\bibitem{Leineweber1999}
Leineweber, D.B.: Efficient reduced {SQP} methods for the optimization of
  chemical processes described by large sparse {DAE} models,
  \emph{Fortschritt-Berichte VDI Reihe 3, Verfahrens\-technik}, vol. 613.
\newblock VDI Verlag, D\"usseldorf
\newblock  (1999)

\bibitem{Leineweber2003}
Leineweber, D.B., Bauer, I., Bock, H.G., Schl\"oder, J.P.: An efficient
  multiple shooting based reduced {SQP} strategy for large-scale dynamic
  process optimization. {P}art {I}: theoretical aspects.
\newblock Computers and Chemical Engineering 27, 157--166
\newblock  (2003)

\bibitem{Li1989}
Li, W.C., Biegler, L.T.: Multistep, {N}ewton-type control strategies for
  constrained nonlinear processes.
\newblock Chem. Eng. Res. Des. 67, 562--577
\newblock  (1989)

\bibitem{Liniger2015}
Liniger, A., Domahidi, A., Morari, M.: Optimization-based autonomous racing of
  1:43 scale {RC} cars.
\newblock Optimal Control Applications and Methods 36(5), 628--647
\newblock  (2015)

\bibitem{Listov2020}
Listov, P., Jones, C.: {PolyMPC}: An efficient and extensible tool for
  real-time nonlinear model predictive tracking and path following for fast
  mechatronic systems.
\newblock Optimal Control Applications and Methods 41(2), 709--727
\newblock  (2020)

\bibitem{Maciejowski2002}
Maciejowski, J.M.: Predictive Control with Constraints.
\newblock Prentice Hall
\newblock  (2002)

\bibitem{MathWorks2005}
MathWorks, T.: The model predictive control toolbox.
\newblock \url{https://mathworks.com/products/mpc.html}
\newblock  (2005)

\bibitem{Mattingley2012}
Mattingley, J., Boyd, S.: {CVXGEN}: A code generator for embedded convex
  optimization.
\newblock Optimization and Engineering pp. 1--27
\newblock  (2012)

\bibitem{Melis2018}
Melis, W., Patrinos, P.: {C} code generation for {NMPC}.
\newblock \url{https://github.com/kul-forbes/nmpc-codegen}
\newblock  (2018)

\bibitem{Nocedal2006}
Nocedal, J., Wright, S.J.: Numerical Optimization, 2 edn.
\newblock Springer Series in Operations Research and Financial Engineering.
  Springer
\newblock  (2006)

\bibitem{Nurkanovic2020a}
Nurkanovi\'c, A., Me{\v s}anovi\'c, A., Zanelli, A., Frey, J., Frison, G.,
  Albrecht, S., Diehl, M.: Real-time nonlinear model predictive control for
  microgrid operation.
\newblock In: Proceedings of the American Control Conference (ACC). Denver, USA
\newblock  (2020).
\newblock (accepted)

\bibitem{Ohtsuka2004}
Ohtsuka, T.: A continuation/{GMRES} method for fast computation of nonlinear
  receding horizon control.
\newblock Automatica 40(4), 563--574
\newblock  (2004)

\bibitem{Pandala2019}
Pandala, A.G., Ding, Y., Park, H.W.: {qpSWIFT}: A real-time sparse quadratic
  program solver for robotic applications.
\newblock {IEEE} Robotics and Automation Letters 4(4), 3355--3362
\newblock  (2019)

\bibitem{Patrinos2014}
Patrinos, P., Bemporad, A.: An accelerated dual gradient-projection algorithm
  for embedded linear model predictive control.
\newblock Automatic Control, IEEE Transactions on 59(1), 18--33
\newblock  (2014)

\bibitem{Qin1996}
Qin, S., Badgwell, T.: {A}n overview of industrial model predictive control
  technology.
\newblock In: J.~Kantor, C.~Garcia, B.~Carnahan (eds.) Fifth International
  Conference on Chemical Process Control -- {CPC V}, pp. 232--256. American
  Institute of Chemical Engineers
\newblock  (1996)

\bibitem{Quirynen2013a}
Quirynen, R., Gros, S., Diehl, M.: Efficient {NMPC} for nonlinear models with
  linear subsystems.
\newblock In: Proceedings of the IEEE Conference on Decision and Control (CDC),
  pp. 5101--5106
\newblock  (2013)

\bibitem{Quirynen2018}
Quirynen, R., Gros, S., Diehl, M.: Inexact {N}ewton-type optimization with
  iterated sensitivities.
\newblock SIAM Journal on Optimization 28(1), 74--95
\newblock  (2018)

\bibitem{Quirynen2017b}
Quirynen, R., Gros, S., Houska, B., Diehl, M.: Lifted collocation integrators
  for direct optimal control in {ACADO} toolkit.
\newblock Mathematical Programming Computation 9(4), 527--571
\newblock  (2017)

\bibitem{Quirynen2017a}
Quirynen, R., Houska, B., Diehl, M.: Efficient symmetric {H}essian propagation
  for direct optimal control.
\newblock Journal of Process Control 50, 19--28
\newblock  (2017)

\bibitem{Quirynen2018a}
Quirynen, R., Knyazev, A., {Di Cairano}, S.: Block structured preconditioning
  within an active-set method for real-time optimal control.
\newblock In: Proceedings of the European Control Conference (ECC)
\newblock  (2018)

\bibitem{Rao1998}
Rao, C.V., Wright, S.J., Rawlings, J.B.: Application of interior-point methods
  to model predictive control.
\newblock Journal of Optimization Theory and Applications 99, 723--757
\newblock  (1998)

\bibitem{Rawlings2017}
Rawlings, J.B., Mayne, D.Q., Diehl, M.M.: Model Predictive Control: Theory,
  Computation, and Design, 2nd edition edn.
\newblock Nob Hill
\newblock  (2017)

\bibitem{Sathya2018}
Sathya, A., Sopasakis, P., Themelis, A., Parys, R.V., Pipeleers, G., Patrinos,
  P.: Embedded nonlinear model predictive control for obstacle avoidance using
  {PANOC}.
\newblock In: Proceedings of the European Control Conference (ECC)
\newblock  (2018)

\bibitem{Schulman2014}
Schulman, J., Duan, Y., Ho, J., Lee, A., Awwal, I., Bradlow, H., Pan, J.,
  Patil, S., Goldberg, K., Abbeel, P.: Motion planning with sequential convex
  optimization and convex collision checking.
\newblock The International Journal of Robotics Research 33(9), 1251--1270
\newblock  (2014)

\bibitem{Shukla2017}
Shukla, H., Khusainov, B., Kerrigan, E., Jones, C.: Software and hardware code
  generation for predictive control using splitting methods.
\newblock In: Proceedings of the IFAC World Congress
\newblock  (2017)

\bibitem{Sopasakis2020}
Sopasakis, P., Fresk, E., Patrinos, P.: {OpEn}: Code generation for embedded
  nonconvex optimization.
\newblock
\newblock In: IFAC World Congress

\bibitem{Steinbach1995}
Steinbach, M.: {F}ast recursive {SQP} methods for large-scale optimal control
  problems.
\newblock {P}h{D} thesis, University of Heidelberg, {IWR}
\newblock  (1995)

\bibitem{Stellato2020}
Stellato, B., Banjac, G., Goulart, P., Bemporad, A., Boyd, S.: {OSQP}: An
  operator splitting solver for quadratic programs.
\newblock Mathematical Programming Computation 12(4), 637--672
\newblock  (2020)

\bibitem{Torrisi2018}
Torrisi, G., Grammatico, S., Smith, R.S., Morari, M.: A projected gradient and
  constraint linearization method for nonlinear model predictive control.
\newblock {SIAM} Journal on Control and Optimization 56(3), 1968--1999
\newblock  (2018)

\bibitem{Ullmann2011}
Ullmann, F.: {FiOrdOs}: A {M}atlab toolbox for {C}-code generation for first
  order methods.
\newblock Master's thesis, ETH Zurich
\newblock  (2011)

\bibitem{Verschueren2016}
Verschueren, R., van Duijkeren, N., Quirynen, R., Diehl, M.: Exploiting
  convexity in direct optimal control: a sequential convex quadratic
  programming method.
\newblock In: Proceedings of the IEEE Conference on Decision and Control (CDC)
\newblock  (2016)

\bibitem{Verschueren2018}
Verschueren, R., Frison, G., Kouzoupis, D., van Duijkeren, N., Zanelli, A.,
  Quirynen, R., Diehl, M.: Towards a modular software package for embedded
  optimization.
\newblock In: Proceedings of the IFAC Conference on Nonlinear Model Predictive
  Control (NMPC)
\newblock  (2018)

\bibitem{Verschueren2017}
Verschueren, R., Zanon, M., Quirynen, R., Diehl, M.: A sparsity preserving
  convexification procedure for indefinite quadratic programs arising in direct
  optimal control.
\newblock SIAM Journal of Optimization 27(3), 2085--2109
\newblock  (2017)

\bibitem{Waechter2006a}
W\"achter, A., Biegler, L.T.: {L}ine {S}earch {F}ilter {M}ethods for
  {N}onlinear {P}rogramming: {M}otivation and {G}lobal {C}onvergence.
\newblock SIAM Journal on Optimization 16, 1--31
\newblock  (2006)

\bibitem{Waechter2006}
W\"achter, A., Biegler, L.T.: On the implementation of an interior-point filter
  line-search algorithm for large-scale nonlinear programming.
\newblock Mathematical Programming 106(1), 25--57
\newblock  (2006)

\bibitem{Wirsching2006}
Wirsching, L., Bock, H.G., Diehl, M.: Fast {NMPC} of a chain of masses
  connected by springs.
\newblock In: Proceedings of the IEEE International Conference on Control
  Applications, Munich, pp. 591--596
\newblock  (2006)

\bibitem{Zanelli2017b}
Zanelli, A., Domahidi, A., Jerez, J.L., Morari, M.: {FORCES NLP}: An efficient
  implementation of interior-point methods for multistage nonlinear nonconvex
  programs.
\newblock International Journal of Control
\newblock  (2017)

\bibitem{Zanelli2019}
Zanelli, A., Quirynen, R., Diehl, M.: Efficient zero-order {NMPC} with
  feasibility and stability guarantees.
\newblock In: Proceedings of the European Control Conference (ECC). Naples,
  Italy
\newblock  (2019)

\bibitem{Zometa2013}
Zometa, P., Kögel, M., Findeisen, R.: {muAO-MPC}: A free code generation tool
  for embedded real-time linear model predictive control.
\newblock In: 2013 American Control Conference, pp. 5320--5325
\newblock  (2013)

\end{thebibliography}


\section*{Appendix A - Case study details}

\section*{Case Study 1: Chain of Masses}

\subsubsection*{System description}

The control objective in this example is to stabilize the motion of a chain of $M=5$ balls with mass $m$ connected by springs to an equilibrium position.
The mass on one end of the chain is fixed at $(0, 0, 0)$.
The mass on the other end can be freely moved.

\begin{figure}
	\centering
	\includegraphics[width=0.7\linewidth]{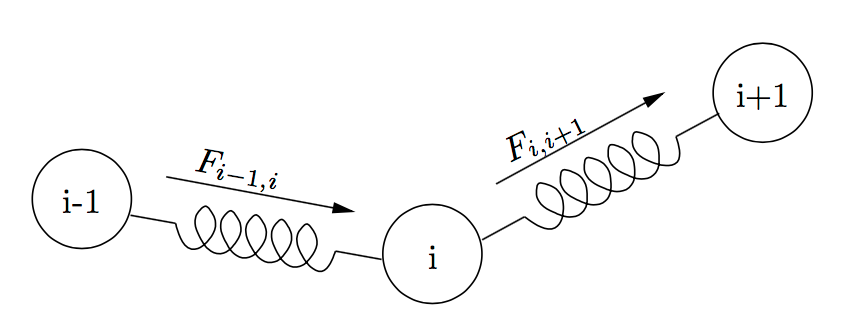}
	\caption{Forces in the springs between the masses in the example of Case Study 1.
		Replicated from \cite{Wirsching2006}.}
	\label{fig:masses}
\end{figure}

Let $p_i$ be the position of mass $i$, for $i=1,\ldots,M$.
The model equations can then be derived as follows.
From Hooke's law, we know that (see Figure~\ref{fig:masses})
\begin{equation*}
	F_{i,i+1} = D\left(1-\frac{L}{\|p_{i+1} - p_{i}\|}\right)(p_{i+1}-p_i),
\end{equation*}
with each spring having spring constant $D$ and rest length $L$.

This allows us to write the equations of motion for the middle balls, which read as
\begin{equation*}
	\ddot{p}_i = \frac{1}{m} (F_{i, i+1} - F_{i-1, i}) + g_z, \quad i=2,\ldots,M-1,
\end{equation*}
with $g_z$ the gravitational acceleration vector.
The control input $ u\in\R^3 $ directly controls the velocity of the free mass:
\begin{equation*}
	\dot{p}_M = u.
\end{equation*}
We now introduce a state space formulation with states
\begin{equation*}
	x^\top = [p_2^\top, p_3^\top, \ldots, p_{M-1}^\top, p_M^\top, v_2^\top, v_3^\top, \ldots, v_{M-1}^\top]\in \R^{n_x}
\end{equation*}
with $n_x=3 \cdot (2\cdot(M-2)+1)$, which results in the following ODE:
\begin{equation}
	\label{eq:chain_of_masses_ode}
	\dot{x} = f(x, u) =  \begin{bmatrix}
		v_2 \\
		\vdots \\ 
		v_{M-1} \\
		u \\
		\frac{1}{m} (F_{2, 3} - F_{1, 2}) + g_z \\
		\vdots \\ 
		\frac{1}{m} (F_{M-1, M} - F_{M-2, M-1}) + g_z \\ 
	\end{bmatrix}.
\end{equation}
We remark that the only nonlinearity is introduced in the calculation of the forces.
The steady state $(x_\mathrm{ss}, u_\mathrm{ss})$ of the system can be found by setting $f(x_\mathrm{ss}, u_\mathrm{ss}) = 0$ for any given $p_{M,\mathrm{ss}}$.
We choose $p_{M,\mathrm{ss}} = [7.5, 0, 0]^\top$ for the experiments.

\subsubsection*{Optimal control problem formulation}

In order to stabilize the motion of the chain of masses to the steady state, we propose the following optimal control problem, obtained by performing a multiple shooting discretization of ODE~\eqref{eq:chain_of_masses_ode}:
\begin{mini}[2]
	{\substack{x_0,\ldots,x_N \\ u_0,\ldots,u_{N-1}}}
	{\sum_{k=0}^{N-1} \begin{bmatrix}x_k-x_\mathrm{ref} \\ u_k-u_\mathrm{ref}\end{bmatrix}^\top \begin{bmatrix} Q & \zeros \\ \zeros & R
		\end{bmatrix} \begin{bmatrix}x_k-x_\mathrm{ref} \\ u_k-u_\mathrm{ref}\end{bmatrix}}
	{}
	{\label{eq:chain_of_masses_ocp}}
	\breakObjective{+ (x_N-x_\mathrm{ref})^\top P (x_N-x_\mathrm{ref})}
	\addConstraint{x_0}{= \overline{x}_0}{}
	\addConstraint{x_{k+1}}{= \phi^x(x_k, u_k), \quad}{k=0,\ldots,N-1}
	\addConstraint{-1}{\leq u_k \leq 1, \quad}{k=0,\ldots,N-1},
\end{mini}
where the initial state $\overline{x}_0$ is the current estimate of the state vector, $\phi: \R^{n_x} \times \R^{3} \rightarrow \R^{n_x}$ is obtained by performing a single RK4 step of length $0.2\,\mathrm{s}$ on ODE~\eqref{eq:chain_of_masses_ode}.
Furthermore, we choose a horizon length $N=40$, the weighting matrices
\begin{align*}
	Q &= \mathrm{diag}(\underbrace{0,\ldots, 0,}_{p_2,\ldots,p_{M-1}} \underbrace{2.5, 2.5, 2.5,}_{p_M} \underbrace{25, \ldots, 25}_{v_2,\ldots,v_{M-1}}), \\
	P &= \mathrm{diag}(\underbrace{0,\ldots, 0,}_{p_2,\ldots,p_{M-1}} \underbrace{10, 10, 10,}_{p_M} \underbrace{0, \ldots, 0}_{v_2,\ldots,v_{M-1}}), \\
	R &= \mathrm{diag}(0.1, 0.1, 0.1),
\end{align*}
the reference control input $ u_\mathrm{ref} = [0, 0, 0]^\top $ and the state reference
\begin{align*}
	x_\mathrm{ref} &= [\underbrace{0,\ldots,0}_{p_2,\ldots,p_{M-1}}, \underbrace{7.5, 0, 0}_{p_M}, \underbrace{0, \ldots, 0}_{v_2,\ldots,v_{M-1}}]^\top.
\end{align*}
The design parameters are chosen as in~\cite{Englert2019} and are summarized in Table~\ref{tab:chain_parameters}.
Note that we did not introduce path constraints or state bounds, since these are not supported by all solvers that we compare to below.

\begin{table}
	\centering
	\caption{Design parameters for the chain of masses case study}
	\label{tab:chain_parameters}
	\begin{tabular}{|c|l|c|}
		\hline 
		Quantity & Description & Value \\ \hline 
		$m$ & mass of one ball & $0.1125\,\mathrm{kg}$ \\ 
		$D$ & spring constant & $0.4\,\mathrm{N/m}$ \\ 
		$L$ & rest length of the springs & $0.1375\,\mathrm{m}$ \\ 
		$g_z$ & gravitational acceleration vector & $[0, 0, -9.81]^\top \, \mathrm{m/s^2}$ \\ 
		$N$ & horizon length & $40$ \\  
		$\Delta t$ & discretization step & $0.2\,\mathrm{s}$ \\ 
		$p_{M,\mathrm{ref}}$ & reference position of free ball & $[7.5, 0, 0]^\top \,\mathrm{m}$ \\
		\hline 
	\end{tabular} 
\end{table}

\section*{Case Study 2: Hessian regularization}

\begin{figure}
	\begin{center}
		\begin{tikzpicture}[scale=0.3,>=stealth]
			\draw[fill=gray!20](-10,0) rectangle (-4,3);
			\draw[fill=black!60] (-9,0) circle(0.5);
			\draw[fill=black!60] (-5,0) circle(0.5);
			\draw[fill=black] (-7,3) circle(0.1);
			\draw[line width=1pt] (-10,10) -- (-7,3);
			\draw[fill=gray!20] (-10,10) circle(1.2);
			\draw[dashed] (-7,1.5) -- (-7,11);
			\draw [->](-7,5.5) arc(90:115:2.1);
			\node at (-8.3,6.3) [right] {$\theta$};
			\node at (-11,1.5){$M$};
			\node at (-10,10){$m$};
			\node at (-9.5,7.2){$l$};
			\draw[->] (-4,1.5) --+(3,0);
			\node at (-2.5,1.5)[above] {$F$};
			\draw[->] (-14,-0.5) -- (-1,-0.5);
			\draw[->] (-14,-0.5) -- (-14,6);
		\end{tikzpicture}
	\end{center}
	\caption{Schematic illustrating the inverted pendulum on top of a cart.}
	\label{fig:cart-pend}
\end{figure}
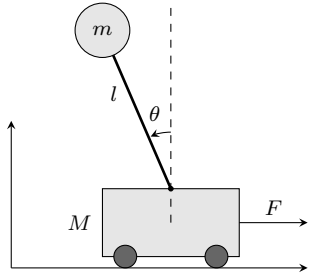

We control a mass on a rod (a pendulum), balanced on a horizontally moving cart, see Figure~\ref{fig:cart-pend}.
The goal is to swing up the pendulum from a stable equilibrium position, namely hanging down.

The dynamics of the cart-pendulum are described by the following ordinary differential equation, where $p, v$ are the horizontal displacement and velocity, respectively, $\theta$ is the angle with the vertical and $\omega$ the corresponding angular velocity:
\begin{subequations}
	\begin{alignat}{2}
		\dot p &= v, \\
		\dot \theta &= \omega, \\
		\dot v &= \frac{-ml\sin(\theta)\dot \theta^2 + mg\cos(\theta)\sin(\theta) + F}{M + m - m(\cos(\theta))^2}, \\
		\dot \omega &= \frac{ -ml\cos(\theta)\sin(\theta)\dot \theta^2 + F\cos(\theta) + (M+m)g\sin(\theta)}{l(M + m - m(\cos(\theta))^2)}.
	\end{alignat}
	\label{eq:explODE_pend}
\end{subequations}

We collect the states in the state vector $x:=[p,\theta,v,\omega]^\top$, the control $u$ is the horizontal force $F$.
Transcribing the continuous-time OCP with multiple shooting gives rise to the following OCP:
\begin{subequations}
	\label{eq:OCP_pendulum}
	\begin{alignat}{3}
		\underset{\substack{x_0, \ldots, x_N, \\ u_0, \ldots, u_{N-1}}}{\mathrm{minimize}} \quad  \sum_{k=0}^{N-1} & \begin{bmatrix}
			x_k \\ u_k
		\end{bmatrix}^\top \begin{bmatrix}
			Q & \zeros \\ \zeros & R
		\end{bmatrix} \begin{bmatrix}
			x_k \\ u_k
		\end{bmatrix} + x_N^\top Q \, x_N\\
		\mathrm{subject~to} \quad x_0 &= \overline{x}_0, \\
		x_{k+1} &= \phi^x_k(x_k,u_k), && k=0,\ldots,N{-}1,\\
		-80 &\leq u_k \leq 80, && k=0,\ldots,N{-}1,
	\end{alignat}
\end{subequations}
where $\phi$ is an RK4 integrator, simulating~\eqref{eq:explODE_pend} over one shooting interval.
The weight matrices are chosen as
\begin{equation*}
	Q=\mathrm{diag}([1000,1000,0.01,0.01]), \quad R=0.01.
\end{equation*}
Because our aim is to swing up the pendulum, we selected strong weights on the position and angle.
The other states and the control are assigned a weak penalty in order to avoid too abrupt swing-ups and to favor smooth trajectories.
Note that the weighting matrices $Q$ and $R$ are tuning parameters used by the control engineer in the design process in order to obtain a desired behavior.
Different choices are therefore equally valid.
The initial value is $\overline{x}_0 = [0,\pi,0,0]^\top$.
We choose $N=100$ shooting intervals of length~$0.01\,\mathrm{s}$.

\section*{Case study 3: Hardware-in-the-loop experiments for an engine control application}

\begin{figure}
	\centering
	\includegraphics[width=0.7\linewidth]{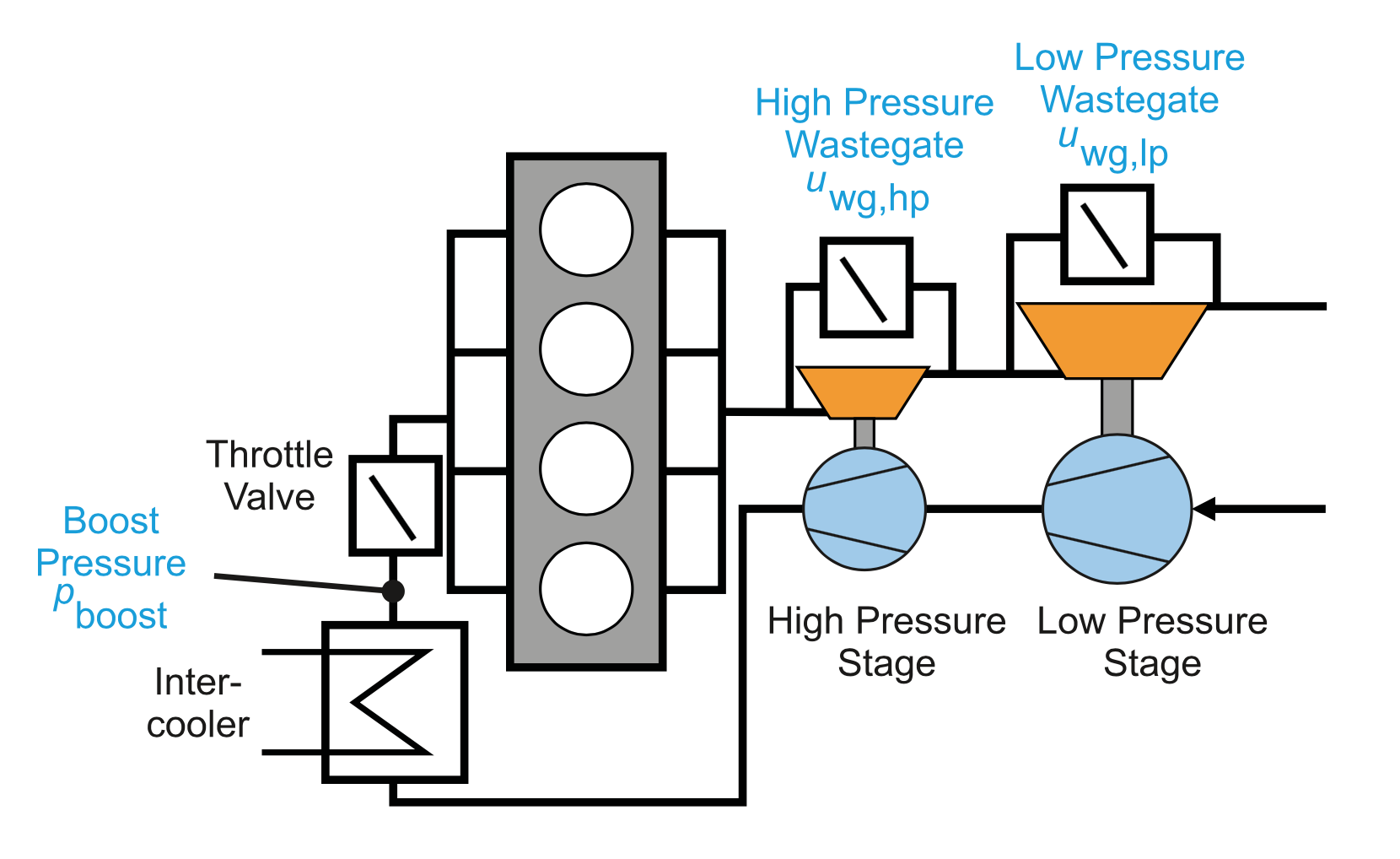}
	\caption{Schematic of the two-stage turbocharging concept.
		Source:~\cite{Albin2017}}
	\label{fig:two-staged-turbocharging}
\end{figure}

Two-stage turbocharging gasoline engines are investigated to overcome the drawbacks of conventional (single-stage) turbocharging.
The main advantage they offer is a better trade-off between short transient times after load changes and a high specific power.
However, the two-stage architecture puts a higher demand on the engine controller, due to the nonlinear nature with cross-couplings in the inputs and the necessity to consider constraints.
NMPC has been proposed as a viable control strategy~\cite{Albin2017}.\\
In Figure~\ref{fig:two-staged-turbocharging}, a sketch of the two-stage turbocharged engine is depicted.
The high-pressure (HP) stage is able to realize fast transients, the low-pressure (LP) stage produces a higher specific power, but with slower dynamics.
The control challenge lies in accurately tracking the boost pressure $p_\mathrm{boost}$, given the highly nonlinear coupling between both stages.\\
For reasons of brevity, we directly present the engine model of~\cite{Albin2017} and refer the interested reader to that work for a derivation.
We model the engine with a set of semi-explicit DAEs.
The differential states consist of $\Pi_{c,\mathrm{lp}}, \Pi_{c,\mathrm{hp}}$ the pressure ratios between input and output of the compressor in the low pressure and high pressure stage, respectively.
The algebraic states are $\Pi_{t,\mathrm{lp}}, \Pi_{t,\mathrm{hp}}$, the pressure ratios on the turbine.
The inputs are the wastegate actuation pulse-width modulated signals $u_\mathrm{wg,lp}, u_\mathrm{wg,hp}$, which take on values between $0\%$ (fully open) and $100\%$ (fully closed).

The resulting DAE system reads as
\begin{align}
	\dot{\Pi}_{c,\mathrm{lp}} &= c_1 (\Pi_{t,\mathrm{lp}}^{1.5} - \Pi_{t,\mathrm{lp}}^{1.25}) \sqrt{\Pi_{t,\mathrm{lp}}^{-1.5} - \Pi_{t,\mathrm{lp}}^{-1.75}} \\ &- c_2 n_\mathrm{eng} \Pi_{c,\mathrm{hp}} (\Pi_{c,\mathrm{lp}}^{1.29} - \Pi_{c,\mathrm{lp}}) \\
	0 &= \Pi_{c,\mathrm{lp}} \Pi_{c,\mathrm{hp}} \\ &- \frac{c_3}{n_\mathrm{eng}} \sqrt{\Pi_{t,\mathrm{lp}}^{0.5} - \Pi_{t,\mathrm{lp}}^{0.25}} \left(\sqrt{\Pi_{t,\mathrm{lp}}} + c_4 \eta(\Pi_{c,\mathrm{lp}} \cdot \Pi_{c,\mathrm{hp}}, u_\mathrm{wg,lp}) \right) \\
	\dot{\Pi}_{c,\mathrm{hp}} &= c_5 (\Pi_{t,\mathrm{hp}}^{1.5} - \Pi_{t,\mathrm{hp}}^{1.25}) \sqrt{\Pi_{t,\mathrm{hp}}^{-1.5} - \Pi_{t,\mathrm{hp}}^{-1.75}} \\ &- c_6 n_\mathrm{eng} \Pi_{c,\mathrm{lp}} (\Pi_{c,\mathrm{hp}}^{1.29} - \Pi_{c,\mathrm{hp}}) \\
	0 &= \Pi_{c,\mathrm{lp}} \Pi_{c,\mathrm{hp}} \\ &- \frac{c_7}{n_\mathrm{eng}} \sqrt{\Pi_{t,\mathrm{hp}}^{0.5} - \Pi_{t,\mathrm{hp}}^{0.25}} \left(\sqrt{\Pi_{t,\mathrm{hp}}} + c_8 (1-u_\mathrm{wg,hp}/100) \right),
	\label{eq:engine_DAE}
\end{align}
with, additionally, $n_\mathrm{eng} = 2000 \, \mathrm{min}^{-1}$ the engine speed.
We model it as a measured disturbance, in this case a constant.
Furthermore, $\eta: \R \times \R \rightarrow \R$ is defined by
\begin{align*}
	\eta(u, v) &= \gamma_1(u) \cdot \gamma_2(v), \\
	\shortintertext{with $\gamma_i: \R \rightarrow \R$:}
	\gamma_i(u) &= b_{1,i} + b_{2,i}\left(1 + e^\frac{-u+b_{3,i}}{b_{4,i}}\right)^{-1}.
\end{align*}
The values of all model parameters can be found in Table~\ref{tab:engine_parameters}.

\begin{table}
	\centering
	\caption{\textsc{Parameter values for the two-stage turbocharged engine model}}
	\label{tab:engine_parameters}
	\begin{tabular}{|c|c|c||c|c|c|}
		\hline 
		Parameter & Unit & Value & Parameter & Unit & Value \\ 
		\hline \hline
		$c_1$ & $\mathrm{-}$ & 25.3 & $b_{1,1}$ & $\mathrm{-}$ & 0 \\ 
		\hline 
		$c_2$ & $\mathrm{min}$ & 0.0034 & $b_{2,1}$ & $\mathrm{-}$ & 1 \\ 
		\hline 
		$c_3$ & $\mathrm{min}^{-1}$ & 7700 & $b_{3,1}$ & $\mathrm{-}$ & 1.49 \\ 
		\hline 
		$c_4$ & $\mathrm{-}$ & 0.6 & $b_{4,1}$ & $\mathrm{-}$ & 0.0377 \\ 
		\hline 
		$c_5$ & $\mathrm{-}$ & 43,6 & $b_{1,2}$ & $\mathrm{-}$ & 67.5 \\ 
		\hline 
		$c_6$ & $\mathrm{min}$ & 0.0092 & $b_{2,2}$ & $\mathrm{-}$ & 4.712 \\ 
		\hline 
		$c_7$ & $\mathrm{min}^{-1}$  & 3600 & $b_{3,2}$ & $\mathrm{-}$ & 1 \\ 
		\hline 
		$c_8$ & $\mathrm{-}$ & 0.9 & $b_{4,2}$ & $\mathrm{-}$ & -1 \\ 
		\hline 
	\end{tabular} 
\end{table}

In order to obtain a smooth control behavior, we include the time derivative of the controls in the optimization formulation, as follows: $\dot{u}_\mathrm{wg,lp} = d_{u, \mathrm{lp}}$, $\dot{u}_\mathrm{wg,hp} = d_{u, \mathrm{hp}}$, and we collect these rates in
\begin{equation*}
	d = \begin{bmatrix}
		d_{u, \mathrm{lp}} \\ d_{u, \mathrm{hp}}
	\end{bmatrix}.
\end{equation*}

We then define the vector of differential states, algebraic states and controls, respectively, as follows:
\begin{equation*}
	x = \begin{bmatrix}
		\Pi_{c,\mathrm{lp}} \\ \Pi_{c,\mathrm{hp}} \\ u_\mathrm{wg,lp} \\ u_\mathrm{wg,hp}
	\end{bmatrix}, \quad z = \begin{bmatrix}
		\Pi_{t,\mathrm{lp}} \\ \Pi_{t,\mathrm{hp}}
	\end{bmatrix}, \quad u = \begin{bmatrix}
		d_{u, \mathrm{lp}} \\ d_{u, \mathrm{hp}}
	\end{bmatrix}.
\end{equation*}

The values of parameters in OCP~\eqref{eq:engine_control} are 
\begin{align*}
	y_{r,k} &= [y_{pr,k}; 1.14; 1.54; 50; 50; 0; 0], \\ W &= \mathrm{diag}([10^3, 10^{-3}, 10^{-3}, 10^{-3}, 10^{-3}, 10^{-4}, 10^{-4}]), \\ 
	y_{r,N} &= [y_{pr,N}; 1.14; 1.54; 50; 50], \\ W_N &= \mathrm{diag}([10^3, 10^{-3}, 10^{-3}, 10^{-3}, 10^{-3}])
\end{align*}
with $y_{pr,k}, k=0,\ldots,N$ varied as in Figure~\ref{fig:enginetracking}.

\end{document}